\documentclass[12pt]{amsart}
\usepackage[margin=2.5cm]{geometry}
\usepackage{tikz-cd}
\usepackage{graphicx}					
\usepackage{amssymb}
\usepackage{amsmath}
\usepackage{relsize}
\usepackage{subfig, graphicx}
\usepackage{mathrsfs}
\usepackage{amssymb, amsthm, indentfirst}
\usepackage{relsize}
\usepackage{hyperref}
\usepackage{stmaryrd}
\usepackage{lineno}
\usepackage{mathabx}

\numberwithin{equation}{section}
\usepackage{color}
\theoremstyle{plain}
\newtheorem{thm}{Theorem}[section]

\newtheorem{conj}[thm]{Conjecture}
\newtheorem{lemma}[thm]{Lemma}

\newtheorem{corollary}[thm]{Corollary}
\theoremstyle{definition}

\newtheorem{rmk}[thm]{Remark}

\newtheorem*{hypothesis*}{Hypothesis}

\author{SHIH-YU CHEN}
\address{Institute of Mathematics~\\Academia Sinica~\\ 6F, Astronomy-Mathematics Building, No.\,1, Sec.\,4, Roosevelt Road, Taipei 10617, Taiwan, ROC}
\email{sychen0626@gate.sinica.edu.tw}

\def\GL{{\rm{GL}}}
\def\GSp{{\rm GSp}}

\def\A{{\mathbb A}}
\def\C{{\mathbb C}}

\def\F{{\mathbb F}}

\def\R{{\mathbb R}}
\def\Q{{\mathbb Q}}
\def\Z{{\mathbb Z}}

\def\<{\langle}
\def\>{\rangle}

\def\bp{\begin{pmatrix}}
\def\ep{\end{pmatrix}}

\def\<{\langle}
\def\>{\rangle}

\def\GL{\operatorname{GL}}

\def\GSp{\operatorname{GSp}}

\def\1{\mathbf{1}}

\def\itPi{\mathit{\Pi}}
\def\itPsi{\mathit{\Psi}}
\def\itSigma{\mathit{\Sigma}}

\makeatletter
\newcommand{\exterior}[1]{\mathop{\mathpalette\exterior@{#1}}}
\newcommand{\exterior@}[2]{%
  \raisebox{\depth}{%
  \fontsize{\sf@size}{0}%
  \m@th
  $\ifx#1\displaystyle\textstyle\else#1\fi\bigwedge$}%
  ^{\mspace{-2mu}#2}%
  \kern-\scriptspace
}
\makeatother

\title{On Deligne's conjecture for symmetric fifth $L$-functions of modular forms}
\begin{document}

\begin{abstract}
We prove Deligne's conjecture for symmetric fifth $L$-functions of elliptic newforms of weight greater than $5$.
As a consequence, we establish period relations between motivic periods associated to an elliptic newform and the Betti--Whittaker periods of its symmetric cube functorial lift to $\GL_4$.
\end{abstract}

\maketitle

\section{Introduction}

In \cite{Deligne1979}, Deligne proposed a conjecture on the algebraicity of values of motivic $L$-functions at critical points in terms of motivic periods. 
A special class of examples are the symmetric power $L$-functions associated to modular forms. 
Given a normalized elliptic newform $f$, we can define the twisted symmetric $n$-th $L$-function $L(s,{\rm Sym}^n(f)\otimes\chi)$ for each integer $n \geq 1$ and Dirichlet character $\chi$. Deligne's conjecture for $L(s,{\rm Sym}^n(f)\otimes\chi)$ at critical points was considered by various authors when $n=1,2,3,4,6$ listed as follows:
\begin{itemize}
\item $n=1$: Shimura \cite{Shimura1976}, \cite{Shimura1977}.
\item $n=2$: Sturm \cite{Sturm1980}, \cite{Sturm1989}.
\item $n=3$: Garrett--Harris \cite{GH1993}, C.- \cite{Chen2021d}.
\item $n=4,6$: Morimoto \cite{Morimoto2021}, C.- \cite{Chen2021}, \cite{Chen2021f}.
\end{itemize}
In this paper, we consider the case when $n=5$.
\subsection{Main results}
Let 
\[
f(\tau) = \sum_{n=1}^\infty a_f(n)q^n \in S_\kappa(N,\omega),\quad q=e^{2\pi\sqrt{-1}\,\tau}
\]
be a normalized elliptic newform of weight $\kappa \geq 2$, level $N$, and nebentypus $\omega$.
For each prime number $p \nmid N$, let $\alpha_p,\beta_p$ be roots of the polynomial
\[
X^2-a_f(p)X+p^{\kappa-1}\omega(p). 
\]
For a Dirichlet character $\chi : (\Z/M\Z)^\times \rightarrow \C$, we define the twisted symmetric fifth $L$-function $L(s,{\rm Sym}^5(f)\otimes\chi)$ by an absolutely convergent Euler product
\[
L(s,{\rm Sym}^5(f)\otimes\chi) = \prod_{p}L_p(s,{\rm Sym}^5(f)\otimes\chi)
\]
for ${\rm Re}(s) > 1+\tfrac{5(\kappa-1)}{2}$.
Here the Euler factors are given by
\[
L_p(s,{\rm Sym}^5(f)\otimes\chi) = \prod_{i=0}^5(1-\alpha_p^{5-i}\beta_p^i\cdot \chi(p)p^{-s})^{-1}
\] 
if $p\nmid NM$.
We omit the definition of the Euler factors at $p \mid NM$ but mentioning that as far as algebraicity concerns, it suffices to consider partial $L$-function.
By the result of Barnet-Lamb--Geraghty--Harris--Taylor \cite[Theorem B]{BLGHT2011} (see also \cite{HSBT2010}), the twisted symmetric fifth $L$-function admits meromorphic continuation to $s\in\C$ and satisfies the expected functional equation relating $L(s,{\rm Sym}^5(f)\otimes\chi)$ to $L(5(\kappa-1)+1-s,{\rm Sym}^5(f^\vee)\otimes\chi^{-1})$. Here $f^\vee \in S_\kappa(N,\omega^{-1})$ is the normalized newform dual to $f$.
The Deligne's periods for ${\rm Sym}^5(f)$ are given by (cf.\,\cite[Proposition 7.7]{Deligne1979})
\[
c^\pm({\rm Sym}^5(f)) = (2\pi\sqrt{-1})^{6-3\kappa}\cdot(\sqrt{-1})^{\kappa}\cdot G(\omega)^6\cdot c^\pm(f)^3\cdot \<f,f\>^3.
\]
Here $G(\omega)$ is the Gauss sum of the primitive Dirichlet character associated to $\omega$, $c^\pm(f) \in \C^\times$ are the periods of $f$ as in \cite{Shimura1977}, and $\<f,f\>$ is the Petersson norm of $f$ defined by
\[
\<f,f\> = {\rm vol}(\Gamma_0(N)\backslash\frak{H})^{-1}\int_{\Gamma_0(N)\backslash\frak{H}}|f(\tau)|^2y^{\kappa-2}\,d\tau.
\]
The set of critical points for ${\rm Sym}^5(f)$ is equal to 
\[
\left\{ m \in \Z \, \left\vert \, 2\kappa-1 \leq m \leq 3\kappa-3\right\}\right..
\]
We have the following special case of Deligne's conjecture \cite{Deligne1979} on the algebraicity of the values of $L(s,{\rm Sym}^5(f) \otimes \chi)$ at critical points.
For $\sigma \in {\rm Aut}(\C)$, let ${}^\sigma\!f \in S_{\kappa}(N,{}^\sigma\!\omega)$ be the normalized newform defined so that $a_{{}^\sigma\!f}(n) = \sigma(a_f(n))$ for $n \geq 1$.
\begin{conj}[Deligne]\label{C:Deligne}
Let $\chi$ be a Dirichlet character and $m \in \Z$ a critical point for ${\rm Sym}^5(f)$. For $\sigma\in{\rm Aut}(\C)$, we have
\begin{align*}
\sigma\left(  \frac{L(m,{\rm Sym}^5(f)\otimes\chi)}{(2\pi\sqrt{-1})^{3m}\cdot G(\chi)^3\cdot c^\pm({\rm Sym}^5(f))} \right) = \frac{L(m,{\rm Sym}^5({}^\sigma\!f)\otimes{}^\sigma\!\chi)}{(2\pi\sqrt{-1})^{3m}\cdot G({}^\sigma\!\chi)^3\cdot c^\pm({\rm Sym}^5({}^\sigma\!f))}.
\end{align*}
Here $\pm = (-1)^m\chi(-1)$.
\end{conj}

Following is the main result of this paper.
\begin{thm}\label{T:main}
If $\kappa \geq 6$, then Conjecture \ref{C:Deligne} holds.
\end{thm}

As a consequence of Theorem \ref{T:main}, we prove the following period relations between the periods associated to $f$ and the Betti--Whittaker periods $p({\rm Sym}^3(f),\pm) \in \C^\times$ (cf.\,\cite{Raghuram2009}) of the symmetric cube lifting ${\rm Sym}^3(f)$.

\begin{thm}\label{T:main 2}
Assume $\kappa \geq 6$. There exists a constant $C_\infty \in \C^\times$, depending only on $\kappa$, such that for all $\sigma \in {\rm Aut}(\C)$, we have
\begin{align*}
\sigma \left( \frac{p({\rm Sym}^3(f),\pm)}{C_\infty\cdot G(\omega)^8\cdot c^\pm(f)^2\cdot \<f,f\>^4} \right) = \frac{p({\rm Sym}^3({}^\sigma\!f),\pm)}{C_\infty\cdot G({}^\sigma\!\omega)^8\cdot c^\pm({}^\sigma\!f)^2\cdot \<{}^\sigma\!f,{}^\sigma\!f\>^4}.
\end{align*}
\end{thm}

\subsection{Notation}
Let $\A$ be the ring of adeles of $\Q$ and $\A_f$ its finite part.
Let $\psi = \bigotimes_v \psi_v : \Q\backslash\A \rightarrow \C$ be the additive character defined so that
\begin{align*}
\psi_{p}(x) & = e^{-2\pi \sqrt{-1}\,x} \mbox{ for }x \in \Z[p^{-1}],\\
\psi_{\infty}(x) & = e^{2\pi \sqrt{-1}\,x} \mbox{ for }x \in \R.
\end{align*}
Denote by $\GL_n$ the general linear group over $\Q$. Let $U_n$ be the standard maximal unipotent subgroup of $\GL_n$ consisting of upper unipotent matrices.
Let $\psi_{U_n} : U_n(\Q)\backslash U_n(\A) \rightarrow \C$ be the additive character defined by
\[
\psi_{U_n}(u) = \psi(u_{12}+u_{23}+\cdots+u_{n-1,n}),\quad u=(u_{ij}) \in U_n(\A).
\]
Let $K_n^\circ$ and $K_n$ be closed subgroups of $\GL_n(\R)$ defined by
\[
K_n^\circ = \R_+\times{\rm SO}(n),\quad K_n = \R_+\times{\rm O}(n).
\]
Here we regard the set $\R_+$ of positive real numbers as the topological connected component of the center of $\GL_n(\R)$.
We denote by $\frak{g}_n$, $\frak{k}_n$, and $\frak{so}(n)$ the Lie algebras of $\GL_n(\R)$, $K_n$, and ${\rm SO}(n)$, respectively.

Let $\sigma \in {\rm Aut}(\C)$.
For a complex representation $\itPi$ of a group $G$ on the space ${V}_\itPi$ of $\itPi$, let ${}^\sigma\!\itPi$ be the representation of $G$ defined
\begin{align}\label{E:sigma action}
{}^\sigma\!\itPi(g) = t \circ \itPi(g) \circ t^{-1},
\end{align}
where $t:{V}_\itPi \rightarrow {V}_\itPi$ is a $\sigma$-linear isomorphism. Note that the isomorphism class of ${}^\sigma\!\itPi$ is independent of the choice of $t$.

\section{Algebraicity of the Rankin--Selberg $L$-functions for $\GL_n \times \GL_{n-1}$}

In this section, we recall the algebraicity of critical values of Rankin--Selberg $L$-functions for $\GL_n(\A) \times \GL_{n-1}(\A)$ in terms of Betti--Whittaker periods following \cite{Raghuram2009}, \cite{Grobner2018b}, \cite{GL2021}, \cite{LLS2021}.

\subsection{Tamely isobaric automorphic representations}\label{SS:tamely isobaric}

Let $P = N_PM_P$ be the standard parabolic subgroup of $\GL_n$ with Levi factor $M_P \simeq \GL_{n_1} \times \cdots \times \GL_{n_k}$. For $1 \leq i \leq k$, let $\itPi_{i}$ be an irreducible cuspidal automorphic representation of $\GL_{n_i}(\A)$. Consider the isobaric automorphic representation of $\GL_n(\A)$:
\[
\itPi = \itPi_1 \boxplus \cdots \boxplus \itPi_k.
\]
Following \cite[\S\,4.3]{LLS2021}, we say $\itPi$ is $tamely$ $isobaric$ if there exists $s \in \R$ such that $\itPi_i \otimes |\mbox{ }|_\A^s$ is unitary for all $1 \leq i \leq k$. In this case, $\itPi$ is fully induced, that is,
\[
\itPi \simeq {\rm Ind}_{P(\A)}^{\GL_n(\A)}(\itPi_1 \otimes \cdots \otimes \itPi_k).
\]
We then realize $\itPi$ in the space $\mathcal{A}(\GL_n)$ of automorphic forms on $\GL_n(\A)$ using Eisenstein series as explained in \cite[\S\,1.4.3]{GL2021}. Moreover, $\itPi$ is globally generic, that is, for non-zero automorphic form $\varphi \in \itPi$, the Whittaker function
\[
W_{\varphi}(g) = \int_{U_n(\Q)\backslash U_n(\A)}\varphi(ug)\overline{\psi_{U_n}(u)}\,du^{\rm Tam},\quad g \in \GL_n(\A)
\]
is non-zero. Here $du^{\rm Tam}$ is the Tamagawa measure on $U_n(\A)$. Let $\mathcal{W}(\itPi_\infty)$ and $\mathcal{W}(\itPi_f)$ be the spaces of Whittaker functions of $\itPi_\infty$ and $\itPi_f = \bigotimes_p \itPi_p$ with respect to $\psi_{U_n,\infty}$ and $\psi_{U_n,f}$, respectively. For $W_\infty \in \mathcal{W}(\itPi_\infty)$ and $W_f \in \mathcal{W}(\itPi_f)$, there exists a unique automorphic form $\varphi \in \itPi$ such that
\[
W_\varphi = W_\infty\cdot W_f.
\]
In this way, we obtain a $((\frak{g}_n,K_n)\times \GL_n(\A_f))$-equivariant isomorphism
\begin{align}\label{E:Whittaker isomorphism}
\mathcal{W}(\itPi_\infty) \otimes \mathcal{W}(\itPi_f) \longrightarrow \itPi.
\end{align}

\subsection{Betti--Whittaker periods for $\GL_n$}\label{SS:Whittaker periods}

Let $\itPi = \itPi_1 \boxplus \cdots \boxplus \itPi_k$ be a tamely isobaric automorphic representation of $\GL_n(\A)$ as in \S\,\ref{SS:tamely isobaric}. 
Let $\rho_P$ be the square-root of the modulus character of $P(\A)$. For $1 \leq i \leq k$, let $\rho_i$ be the restriction of $\rho_P$ to the factor $\GL_{n_i}$ of the Levi component $M_P \simeq \GL_{n_1} \times \cdots \times \GL_{n_k}$.
Assume further that $\itPi$ is $cohomological$, that is, there exists an irreducible algebraic representation $M$ of $\GL_n$ such that the $(\frak{g}_n,K_n^\circ)$-cohomology 
\[
H^\bullet(\frak{g}_n,K_n^\circ;\itPi_\infty\otimes M_\C) \neq 0.
\]
Here $M_\C = M\otimes_\Q\C$. Note that $M$ is uniquely determined by $\itPi_\infty$ and $\bullet = b_n = \lfloor \tfrac{n^2}{4} \rfloor$ is the least integer such that the $(\frak{g}_n,K_n^\circ)$-cohomology is non-vanishing (cf.\,\cite[\S\,3]{Clozel1990}).
Moreover, we have
\[
{\rm dim}_\C \,H^{b_n}(\frak{g}_n,K_n^\circ;\itPi_\infty\otimes M_\C) = \begin{cases}
1 & \mbox{ if $n$ is odd},\\
2 & \mbox{ if $n$ is even}.
\end{cases}
\]
When $n$ is odd, we fix a generator 
\begin{align}\label{generator 1}
[\itPi_\infty] \in H^{b_n}(\frak{g}_n,K_n^\circ;\itPi_\infty\otimes M_\C).
\end{align}
When $n$ is even, the $\pm1$-eigenspace $H^{b_n}(\frak{g}_n,K_n^\circ;\itPi_\infty\otimes M_\C)[\pm]$ under the action of $K_n / K_n^\circ$ is $1$-dimensional. We fix a generator 
\begin{align}\label{generator 2}
[\itPi_\infty]_\pm \in H^{b_n}(\frak{g}_n,K_n^\circ;\itPi_\infty\otimes M_\C)[\pm].
\end{align}

\subsubsection{Rational structure via Whittaker model}
Let $\sigma \in {\rm Aut}(\C)$. For $1 \leq i \leq k$, let ${}^\sigma(\itPi_i \otimes \rho_i^{-1})$ be the irreducible admissible $((\frak{g}_n,K_n) \times \GL_n(\A_f))$-module defined by
\[
{}^\sigma(\itPi_i \otimes \rho_i^{-1}) = (\itPi_{i,\infty}\otimes \rho_{i,\infty}^{-1}) \otimes {}^\sigma(\itPi_{i,f} \otimes \rho_{i,f}^{-1}).
\]
By \cite[III, Theorem 3.3]{BW2000}, $\itPi_i \otimes \rho_i^{-1}$ is coholomogical. In particular, ${}^\sigma(\itPi_i \otimes \rho_i^{-1})$ is cohomological and cuspidal by the result of Clozel \cite[Th\'eor\`eme 3.19]{Clozel1990}. Let ${}^\sigma\itPi$ be the isobaric automorphic representation defined by
\[
{}^\sigma\!\itPi = \bigboxplus_{i=1}^k {}^\sigma(\itPi_i \otimes \rho_i^{-1})\otimes \rho_i.
\]
It is clear that ${}^\sigma\!\itPi$ is cohomological and tamely isobaric. Moreover, we have $({}^\sigma\!\itPi)_f = {}^\sigma(\itPi_f)$ (cf.\,\cite[Lemma 1.2]{Grobner2018b}).
We write ${}^\sigma\!\itPi_f = ({}^\sigma\!\itPi)_f$. 
Let $\Q(\itPi)$ be the rationality field of $\itPi$, which is the fixed field of $\{\sigma \in {\rm Aut}(\C)\,\vert\,{}^\sigma\!\itPi = \itPi\}$.
Note that $\Q(\itPi)$ is equal to the composite of the rationality fields of $\Q(\itPi_i)$ for $1 \leq i \leq k$. In particular, $\Q(\itPi)$ is a number field.
Let $t_\sigma : \mathcal{W}(\itPi_f) \rightarrow \mathcal{W}({}^\sigma\!\itPi_f)$ be the $\sigma$-linear $\GL_n(\A_f)$-equivariant isomorphism defined by
\[
t_\sigma W(g) = \sigma\left(W({\rm diag}(u_\sigma^{-n+1},u_\sigma^{-n+2},\cdots,1)g)\right),\quad g \in \GL_n(\A_f).
\]
Here $u_\sigma \in \widehat{\Z}^\times$ is the unique element such that $\sigma(\psi(x)) = \psi(u_\sigma x)$ for all $x \in \A_f$.
We thus obtain a $\Q(\itPi)$-rational structure on $\mathcal{W}(\itPi_f)$ given by taking the Galois invariants:
\begin{align}\label{E:rational structure 1}
\mathcal{W}(\itPi_f)^{{\rm Aut}(\C/\Q(\itPi))} = \left.\left\{W\in \mathcal{W}(\itPi_f)\,\right\vert\,t_\sigma W=M\mbox{ for $\sigma \in {\rm Aut}(\C/\Q(\itPi))$}\right\}.
\end{align}

\subsubsection{Rational structure via sheaf cohomology}

Consider the orbifold
\[
\mathcal{X}_n = \GL_n(\Q)\backslash \GL_n(\A) /K_n^\circ. 
\]
The algebraic representation $M_\C$ defines a locally constant sheaf
$\mathcal{M}_{\C}$ of $\C$-vector spaces on $\mathcal{X}_n$ (cf.\,\cite[\S\,2.2.7]{HR2020}).
We denote by
\[
H^\bullet(\mathcal{X}_n,\mathcal{M}_\C)
\]
the sheaf cohomology of $\mathcal{M}_\C$.
For $\sigma \in {\rm Aut}(\C)$, the canonical $\sigma$-linear isomorphism $M_\C \rightarrow M_\C$ naturally induces a $\sigma$-linear isomorphism
\[
\sigma^* : H^\bullet(\mathcal{X}_n,\mathcal{M}_\C) \longrightarrow H^\bullet(\mathcal{X}_n,\mathcal{M}_\C).
\]
Assume $n$ is odd. Following the proof of \cite[Propositions 1.6 and 1.7]{Grobner2018b} (cf.\,\cite[Theorem 7.23]{GR2014b} and \cite[Proposition 4.3]{LLS2021}), there exists a $\GL_n(\A_f)$-equivariant injection
\[
\Psi_\itPi : H^{b_n}(\frak{g}_n,K_n^\circ;\itPi \otimes M_\C) \longrightarrow H^{b_n}(\mathcal{X}_n,\mathcal{M}_\C)
\]
such that the image of $\sigma^*\circ \Psi_\itPi$ is the image of $\Psi_{{}^\sigma\!\itPi}$. This induces a $\sigma$-linear $\GL_n(\A_f)$-equivariant isomorphism
\[
\sigma^* : H^{b_n}(\frak{g}_n,K_n^\circ;\itPi \otimes M_\C) \longrightarrow H^{b_n}(\frak{g}_n,K_n^\circ;{}^\sigma\!\itPi \otimes M_\C).
\]
We thus obtain a $\Q(\itPi)$-rational structure on $H^{b_n}(\frak{g}_n,K_n^\circ;\itPi \otimes M_\C)$ given by taking the Galois invariants:
\begin{align}\label{E:rational structure 2}
\begin{split}
&H^{b_n}(\frak{g}_n,K_n^\circ;\itPi \otimes M_\C)^{{\rm Aut}(\C/\Q(\itPi))}\\
& = \left.\left\{c\in H^{b_n}(\frak{g}_n,K_n^\circ;\itPi \otimes M_\C)\,\right\vert\,\sigma^* c=c\mbox{ for $\sigma \in {\rm Aut}(\C/\Q(\itPi))$}\right\}.
\end{split}
\end{align}
Let 
\[
\Phi_\itPi : \mathcal{W}(\itPi_f) \longrightarrow H^{b_n}(\frak{g}_n,K_n^\circ;\itPi \otimes M_\C)
\]
be the $\GL_n(\A_f)$-equivariant isomorphism defined as follows: For $W \in \mathcal{W}(\itPi_f)$, we have
\[
[\itPi_\infty]\otimes W \in H^{b_n}(\frak{g}_n,K_n^\circ; \mathcal{W}(\itPi_\infty) \otimes \mathcal{W}(\itPi_f)\otimes M_\C).
\] 
Then $\Phi_\itPi(W)$ is the image of $[\itPi_\infty]\otimes W$ under the $\GL_n(\A_f)$-equivariant isomorphism induced by the isomorphism (\ref{E:Whittaker isomorphism}).
When $n$ is even, similarly we have a $\Q(\itPi)$-rational structure $H^{b_n}(\frak{g}_n,K_n^\circ;\itPi \otimes M_\C)[\pm]^{{\rm Aut}(\C/\Q(\itPi))}$ on $H^{b_n}(\frak{g}_n,K_n^\circ;\itPi \otimes M_\C)[\pm]$ and a $\GL_n(\A_f)$-equivariant isomorphism
\[
\Phi_{\itPi,\pm} : \mathcal{W}(\itPi_f) \longrightarrow H^{b_n}(\frak{g}_n,K_n^\circ;\itPi \otimes M_\C)[\pm].
\]
Comparing the $\Q(\itPi)$-rational structures given by (\ref{E:rational structure 1}) and (\ref{E:rational structure 2}), we have the following lemma/definition of the Betti--Whittaker periods of $\itPi$.

\begin{lemma}
There exists $p(\itPi) \in \C^\times$ if $n$ is odd and $p(\itPi,\pm)\in \C^\times$ if $n$ is even, unique up to $\Q(\itPi)^\times$, such that
\[
\frac{\Phi_\itPi(\mathcal{W}(\itPi_f)^{{\rm Aut}(\C/\Q(\itPi))})}{p(\itPi)} = H^{b_n}(\frak{g}_n,K_n^\circ;\itPi \otimes M_\C)^{{\rm Aut}(\C/\Q(\itPi))}
\]
if $n$ is odd and
\[
\frac{\Phi_{\itPi,\pm}(\mathcal{W}(\itPi_f)^{{\rm Aut}(\C/\Q(\itPi))})}{p(\itPi,\pm)} = H^{b_n}(\frak{g}_n,K_n^\circ;\itPi \otimes M_\C)[\pm]^{{\rm Aut}(\C/\Q(\itPi))}
\]
if $n$ is even.
Moreover, we can normalize the periods so that
\[
\sigma^* \left(\frac{\Phi_{\itPi}(W)}{p(\itPi)}\right)  = \frac{\Phi_{{}^\sigma\!\itPi}(t_\sigma W)}{p({}^\sigma\!\itPi)}
\]
if $n$ is odd and
\[
\sigma^* \left(\frac{\Phi_{\itPi,\pm}(W)}{p(\itPi,\pm)}\right) = \frac{\Phi_{{}^\sigma\!\itPi,\pm}(t_\sigma W)}{p({}^\sigma\!\itPi,\pm)}
\]
if $n$ is even for all $W \in \mathcal{W}(\itPi_f)$ and $\sigma \in {\rm Aut}(\C)$.
\end{lemma}

\subsection{Algebraicity for $\GL_n \times \GL_{n-1}$}

For isobaric automorphic representations $\itPi_n$ and $\itPi_m$ of $\GL_n(\A)$ and $\GL_m(\A)$, respectively, let
\[
L(s,\itPi_n \times \itPi_m)
\]
be the Rankin--Selberg $L$-function for $\itPi_n \times \itPi_m$. We denote by $L^{(\infty)}(s,\itPi_n \times \itPi_m)$ the $L$-function obtained by excluding the archimedean $L$-factor.

Let $\itPi_{n,\infty}$ and $\itPi_{n-1,\infty}$ be irreducible admissible generic essentially unitary $(\frak{g}_{n},K_n^\circ)$ and $(\frak{g}_{n-1},K_{n-1}^\circ)$ modules, respectively. Let $M_n$ and $M_{n-1}$ be the irreducible algebraic representations of $\GL_{n}$ and $\GL_{n-1}$ such that
\[
H^{b_n}(\frak{g}_{n},K_n^\circ;\itPi_{n,\infty}\otimes M_{n,\C})\neq 0,\quad H^{b_{n-1}}(\frak{g}_{n-1},K_{n-1}^\circ;\itPi_{n-1,\infty}\otimes M_{n-1,\C})\neq 0.
\]
Here $M_{n,\C} = M_n \otimes_\Q \C$ and $M_{n-1,\C} = M_{n-1} \otimes_\Q \C$. Note that $M_n$ and $M_{n-1}$ are uniquely determined by $\itPi_{n,\infty}$ and $\itPi_{n-1,\infty}$ (cf.\,\cite[\S\,3]{Clozel1990}). We say $(\itPi_{n,\infty},\itPi_{n-1,\infty})$ is $balanced$ if there exists $m \in \Z$ such that
\begin{align}\label{E:balanced}
{\rm Hom}_{\GL_{n-1}}(M_n \otimes M_{n-1},{\det}^m) \neq 0.
\end{align}
In this case, $m' \in \Z$ satisfying (\ref{E:balanced}) if and only if $L(s,\itPi_{n,\infty}\times \itPi_{n-1,\infty})$ and $L(1-s,\itPi_{n,\infty}^\vee\times \itPi_{n-1,\infty}^\vee)$ are holomorphic at $s= m' + \tfrac{1}{2}$ (cf.\,\cite[Theorem 2.21]{Raghuram2016}).
Fix $m \in \Z$ satisfying (\ref{E:balanced}). We define a quantity $Z(m,\itPi_{n,\infty},\itPi_{n-1,\infty}) \in \C$ as follows: Let $\iota : \GL_{n-1} \rightarrow \GL_n$ be the embedding defined by 
\[
\iota(g) = \bp g & 0 \\ 0 & 1\ep.
\]
\begin{itemize}
\item Fix a non-zero $\xi_m \in {\rm Hom}_{\GL_{n-1}}(M_n \otimes M_{n-1},{\det}^m)$.
\item For $1 \leq i ,j \leq n-1$, let $e_{ij} \in \frak{g}_{n-1}$ be the $n-1$ by $n-1$ matrix with $1$ in the $(i,j)$-entry and zeros otherwise.
By abuse of notation, the image of $e_{ij} \in \frak{g}_{n-1}$ in $\frak{g}_{n-1} / \frak{so}({n-1})$ is denoted by the same symbol $e_{ij}$.
Then the generator
\[
(\wedge_{i=1}^{n-1} \,e_{ii}^*) \wedge (\wedge_{1\leq i <j \leq n-1}\,e_{ij}^* ) \in \exterior{n(n-1)/2}(\frak{g}_{n-1,\C}/\frak{so}(n-1)_\C)^*
\]
defines an invariant measure on $\GL_{n-1}(\R)/{\rm SO}(n-1)$. By pull-back, this then defines a Haar measure on $\GL_{n-1}(\R)$ by requiring ${\rm vol}({\rm SO}(n-1))=1$.
For $X \in \exterior{b_n}(\frak{g}_{n,\C}/\frak{k}_{n,\C})^*$ and $Y \in \exterior{b_{n-1}}(\frak{g}_{n-1,\C}/\frak{k}_{n-1,\C})^*$, let ${\sf s}(X,Y) \in \C$ defined so that
\[
\iota^*(X) \wedge {\rm pr}(Y) = {\sf s}(X,Y)\cdot (\wedge_{i=1}^{n-1} \,e_{ii}^*) \wedge (\wedge_{1\leq i <j \leq n-1}\,e_{ij}^* ).
\]
Here ${\rm pr}:\exterior{b_{n-1}}(\frak{g}_{n-1,\C}/\frak{k}_{n-1,\C})^* \rightarrow \exterior{b_{n-1}}(\frak{g}_{n-1,\C}/\frak{so}(n-1)_\C)^*$ is the natural surjection and $\iota^* : \exterior{b_n}(\frak{g}_{n,\C}/\frak{k}_{n,\C})^* \rightarrow \exterior{b_{n}}(\frak{g}_{n-1,\C}/\frak{so}(n-1)_\C)^*$ is induced by the embedding $\iota$.

\item For $W_n\in \mathcal{W}(\itPi_{n,\infty})$ and $W_{n-1}\in \mathcal{W}(\itPi_{n-1,\infty})$, we define local zeta integral 
\[
Z(s,W_n,W_{n-1}) = \int_{U_{n-1}(\R)\backslash\GL_{n-1}(\R)}W_n(\iota(g))W_{n-1}(g)|\det g|^{s-1/2}\,d\overline{g}.
\]
Here $d\overline{g}$ is the quotient measure induced by the Lebesgue measure on $U_{n-1}(\R) \simeq \R^{(n-1)(n-2)/2}$ and the Haar measure on $\GL_{n-1}(\R)$ described above.
The integral converges absolutely for ${\rm Re}(s)$ sufficiently large and admits meromorphic continuation to $s \in \C$. Moreover, the ratio
\[
\frac{Z(s,W_n,W_{n-1})}{L(s,\itPi_{n,\infty}\times \itPi_{n-1,\infty})}
\]
is entire (cf.\,\cite{JS1990}, \cite{Jacquet2009}). In particular, $Z(s,W_n,W_{n-1})$ is holomorphic at $s = m+\tfrac{1}{2}$.
\end{itemize}
We then define the $K_{n-1}^\circ$-equivariant bilinear pairing $\<\,,\,\>_m $ on
\[
\left(\mathcal{W}(\itPi_{n,\infty})\otimes \exterior{b_n}(\frak{g}_{n,\C}/\frak{k}_{n,\C})^* \otimes M_{n,\C}\right) \times \left(\mathcal{W}(\itPi_{n-1,\infty})\otimes \exterior{b_{n-1}}(\frak{g}_{n-1,\C}/\frak{k}_{n-1,\C})^* \otimes M_{n-1,\C}\right)
\]
by
\[
\left\<W_n \otimes X \otimes {\bf v},\, W_{n-1}\otimes Y \otimes {\bf v}'\right\>_m = Z(m+\tfrac{1}{2},W_n,W_{n-1})\cdot {\sf s}(X,Y)\cdot \xi_m({\bf v}\otimes{\bf v}').
\]
By restriction, we obtain a bilinear pairing $\<\,,\,\>_m$ on
\[
H^{b_n}(\frak{g}_{n},K_{n}^\circ;\itPi_{n,\infty}\otimes M_{n,\C})\times H^{b_{n-1}}(\frak{g}_{n-1},K_{n-1}^\circ;\itPi_{n-1,\infty}\otimes M_{n-1,\C}).
\]
Let ${\sf w}(\itPi_{n,\infty})$ and ${\sf w}(\itPi_{n-1,\infty})$ be the integers such that 
\[
|\omega_{\itPi_{n,\infty}}| = |\mbox{ }|^{n{\sf w}(\itPi_{n,\infty})/2},\quad |\omega_{\itPi_{n-1,\infty}}| = |\mbox{ }|^{(n-1){\sf w}(\itPi_{n-1,\infty})/2}.
\]
When $n$ is odd (resp.\,even), ${\sf w}(\itPi_{n,\infty})$ (resp.\,${\sf w}(\itPi_{n-1,\infty})$) must be even and we put
\[
\varepsilon(\itPi_{n,\infty}) = (-1)^{{\sf w}(\itPi_{n,\infty})/2}\omega_{\itPi_{n,\infty}}(-1) \quad ({\rm resp.}\,\varepsilon(\itPi_{n-1,\infty}) = (-1)^{{\sf w}(\itPi_{n-1,\infty})/2}\omega_{\itPi_{n-1,\infty}}(-1)).
\]
We fix generators in the $(\frak{g}_n,K_n^\circ)$ and $(\frak{g}_{n-1},K_{n-1}^\circ)$ cohomologies as in (\ref{generator 1} and \ref{generator 2}). Define $Z(m,\itPi_{n,\infty},\itPi_{n-1,\infty}) \in \C$ by
\[
Z(m,\itPi_{n,\infty},\itPi_{n-1,\infty}) = \begin{cases}
\left\<[\itPi_{n,\infty}],\,[\itPi_{n-1,\infty}]_{(-1)^{1+m}\varepsilon(\itPi_{n,\infty})}\right\>_m & \mbox{ if $n$ is odd},\\
\left\<[\itPi_{n,\infty}]_{(-1)^m\varepsilon(\itPi_{n-1,\infty})},\,[\itPi_{n-1,\infty}]\right\>_m & \mbox{ if $n$ is even}.
\end{cases}
\]
We have the following non-vanishing result due to Sun \cite{Sun2017}.
\begin{thm}[Sun]
Assume $(\itPi_{n,\infty},\itPi_{n-1,\infty})$ is balanced and $m\in\Z$ satisfying (\ref{E:balanced}).
We have 
\[
Z(m,\itPi_{n,\infty},\itPi_{n-1,\infty}) \neq 0.
\]
\end{thm}

\begin{rmk}
In \cite{LLS2021}, Li--Liu--Sun establish period relations for the ratios $\frac{Z(m,\itPi_{n,\infty},\itPi_{n-1,\infty})}{Z(m',\itPi_{n,\infty},\itPi_{n-1,\infty})}$.
\end{rmk}

In the following theorem, we recall a generalization of the result of Raghuram \cite{Raghuram2009} on the algebraicity of Rankin--Selberg $L$-functions for $\GL_n \times \GL_{n-1}$ in terms of Betti--Whittaker periods and the above archimedean quantity.
\begin{thm}\label{T:Raghuram}
Let $\itPi_n$ (resp.\,$\itPi_{n-1}$) be cohomological irreducible cuspidal (resp.\,tamely isobaric) automorphic representation of $\GL_n(\A)$ (resp.\,$\GL_{n-1}(\A)$) such that $(\itPi_{n,\infty},\itPi_{n-1,\infty})$ is balanced.
For a critical point $m+\tfrac{1}{2} \in \Z+\tfrac{1}{2}$ of $L(s,\itPi_n \times \itPi_{n-1})$, we put
\begin{align*}
&p(m,\itPi_n \times \itPi_{n-1})\\
& = Z(m,\itPi_{n,\infty},\itPi_{n-1,\infty})^{-1}\cdot G(\omega_{\itPi_{n-1}})\cdot\begin{cases}
p(\itPi_n)p(\itPi_{n-1},(-1)^{1+m}\varepsilon(\itPi_{n,\infty})) & \mbox{ if $n$ is odd},\\
p(\itPi_n,(-1)^m\varepsilon(\itPi_{n-1,\infty}))p(\itPi_{n-1}) & \mbox{ if $n$ is even}.
\end{cases}
\end{align*}
Then 
\[
\sigma \left( \frac{L^{(\infty)}(m+\tfrac{1}{2},\itPi_n \times \itPi_{n-1})}{p(m,\itPi_n \times \itPi_{n-1})} \right) = \frac{L^{(\infty)}(m+\tfrac{1}{2},{}^\sigma\!\itPi_n \times {}^\sigma\!\itPi_{n-1})}{p(m,{}^\sigma\!\itPi_n \times {}^\sigma\!\itPi_{n-1})}
\]
for all $\sigma \in {\rm Aut}(\C)$.
\end{thm}

\begin{rmk}
In \cite{Raghuram2009}, the result is proved under the assumption that $\itPi_n$ and $\itPi_{n-1}$ are both cuspidal. For tamely isobaric $\itPi_{n-1}$, the proof goes without change, except we consider the image of $H^{b_{n-1}}(\frak{g}_n,K_n^\circ; \itPi_{n-1}\otimes M_{n-1,\C})$ in $H^{b_{n-1}}(\mathcal{X}_{n-1},\mathcal{M}_{n-1,\C})$ instead of the cohomology with compact support when $\itPi_{n-1}$ is cuspidal. We refer to \cite[\S\,6]{LLS2021} for general base field (see also \cite{GH2016}, \cite{Grobner2018b} when the base field is a CM-field).
\end{rmk}

As an immediate consequence of the theorem, we have the following corollary on the algebraicity of ratios of product of critical values of the Rankin--Selberg $L$-functions.

\begin{corollary}\label{C:ratio of L-values}
Let $\itPi_n$ and $\itPi_n'$ (resp.\,$\itPi_{n-1}$ and $\itPi_{n-1}'$) be cohomological irreducible cuspidal (resp.\,tamely isobaric) automorphic representations of $\GL_n(\A)$ (resp.\,$\GL_{n-1}(\A)$) satisfying the following conditions:
\begin{itemize}
\item $(\itPi_{n,\infty},\itPi_{n-1,\infty})$ is balanced.
\item $\itPi_{n,\infty} = \itPi_{n,\infty}'$ and $\itPi_{n-1,\infty} = \itPi_{n-1,\infty}'$.
\end{itemize}
For $\sigma \in {\rm Aut}(\C)$ and $m+\tfrac{1}{2}$ a critical point such that $L(m+\tfrac{1}{2},\itPi_n \times \itPi_{n-1})\cdot L(m+\tfrac{1}{2},\itPi_n' \times \itPi_{n-1}')$ is non-zero,
we have
\begin{align*}
&\sigma \left( \frac{L^{(\infty)}(m+\tfrac{1}{2},\itPi_n \times \itPi_{n-1}')\cdot L^{(\infty)}(m+\tfrac{1}{2},\itPi_n' \times \itPi_{n-1})}{L^{(\infty)}(m+\tfrac{1}{2},\itPi_n \times \itPi_{n-1})\cdot L^{(\infty)}(m+\tfrac{1}{2},\itPi_n' \times \itPi_{n-1}')} \right)\\
& = \frac{L^{(\infty)}(m+\tfrac{1}{2},{}^\sigma\!\itPi_n \times {}^\sigma\!\itPi_{n-1}')\cdot L^{(\infty)}(m+\tfrac{1}{2},{}^\sigma\!\itPi_n' \times {}^\sigma\!\itPi_{n-1})}{L^{(\infty)}(m+\tfrac{1}{2},{}^\sigma\!\itPi_n \times {}^\sigma\!\itPi_{n-1})\cdot L^{(\infty)}(m+\tfrac{1}{2},{}^\sigma\!\itPi_n' \times {}^\sigma\!\itPi_{n-1}')}.
\end{align*}
\end{corollary}

\subsection{Case $n=2$}
Let $\itPi$ be a cohomological irreducible cuspidal automorphic representation of $\GL_2(\A)$ with central character $\omega_\itPi$. 
Then there exist $\kappa \geq 2$ and ${\sf w} \in \Z$ such that $\kappa \equiv {\sf w}\,({\rm mod}\,2)$ and 
\[
\itPi_\infty = D_\kappa \otimes |\mbox{ }|^{{\sf w}/2}.
\]
Let $M_\mu$ be the irreducible algebraic representation of $\GL_2$ of highest weight
\[
\mu = \left( \tfrac{\kappa-2-{\sf w}}{2},\tfrac{-\kappa+2-{\sf w}}{2}\right).
\]
Following the normalization in \cite[(3.6) and (3.7)]{RT2011}, we define the generators
\[
[\itPi_\infty]_\pm\in  H^1(\frak{g}_2,K_2^\circ;\mathcal{W}(\itPi_\infty) \otimes M_{\mu,\C})[\pm].
\]
In the following theorem, we recall a refinement of Theorem \ref{T:Raghuram} when $n=2$ due to Raghuram--Tanabe \cite{RT2011} and a period relation.
Let $f_\itPi \in \itPi$ be the normalized newform of $\itPi$ in the sense of Casselman \cite{Casselman1973}.
Let $\Vert f_\itPi \Vert$ be the Petersson norm of $f_\itPi$ defined by
\begin{align}\label{E:Petersson norm}
\Vert f_\itPi \Vert = \int_{\A^\times\GL_2(\Q)\backslash\GL_2(\A)}|f_\itPi(g)|^2|\det g|_\A^{-{\sf w}}\,dg^{\rm Tam}.
\end{align}
Here $dg^{\rm Tam}$ is the Tamagawa measure on $\A^\times\backslash \GL_2(\A)$.

\begin{thm}\label{T:GL_2}
\noindent
\begin{itemize}
\item[(1)]{\rm (Raghuram--Tanabe)} Let $\chi$ be a finite order Hecke character of $\A^\times$ and $m+\tfrac{1}{2}$ a critical point of $L(s,\itPi\otimes\chi)$. For $\sigma \in {\rm Aut}(\C)$, we have
\begin{align*}
&\sigma\left( \frac{L^{(\infty)}(m+\tfrac{1}{2},\itPi \otimes \chi)}{(2\pi\sqrt{-1})^{m+(\kappa+{\sf w})/2}\cdot G(\chi) \cdot p(\itPi,(-1)^m\cdot\chi_\infty(-1))}\right)\\
&= \frac{L^{(\infty)}(m+\tfrac{1}{2},{}^\sigma\!\itPi \otimes {}^\sigma\!\chi)}{(2\pi\sqrt{-1})^{m+(\kappa+{\sf w})/2}\cdot G({}^\sigma\!\chi) \cdot p({}^\sigma\!\itPi,(-1)^m\cdot\chi_\infty(-1))}.
\end{align*}
\item[(2)] For $\sigma \in {\rm Aut}(\C)$, we have
\begin{align*}
\sigma\left(\frac{p(\itPi,+)\cdot p(\itPi,-)}{(2\pi\sqrt{-1})\cdot(\sqrt{-1})^{{\sf w}}\cdot G(\omega_\itPi)\cdot\Vert f_\itPi \Vert}\right) = \frac{p({}^\sigma\!\itPi,+)\cdot p({}^\sigma\!\itPi,-)}{(2\pi\sqrt{-1})\cdot(\sqrt{-1})^{{\sf w}}\cdot G({}^\sigma\!\omega_\itPi)\cdot\Vert f_{{}^\sigma\!\itPi} \Vert}.
\end{align*}
\end{itemize}
\end{thm}

\begin{proof}
We refer to \cite[Lemma 4.1]{Chen2020} for the proof of the second assertion.
\end{proof}

\section{Proof of main results}

\subsection{Some algebraicity results}

Firstly we recall some algebraicity results in the literature which will be used in the proof of our main result.

\begin{thm}\label{T:algebraicity results}
For $i=1,2,3$, let $\itPi_i$ be a cohomological irreducible cuspidal automorphic representation of $\GL_2(\A)$ with $\itPi_{i,\infty} = D_{\kappa_i}\otimes|\mbox{ }|^{{\sf w}_i/2}$ for some $\kappa_i \equiv {\sf w}_i \,({\rm mod}\,2)$.
\begin{itemize}
\item[(1)]{\rm (Shimura\,\cite[Theorem 3]{Shimura1976})} Assume $\kappa_1 >\kappa_2$. For $\sigma \in {\rm Aut}(\C)$ and $m \in \Z$ critical for $L(s,\itPi_1\times\itPi_2)$,
we have
\begin{align*}
&\sigma \left( \frac{L^{(\infty)}(m,\itPi_1 \times \itPi_2)}{(2\pi\sqrt{-1})^{2m+\kappa_1+{\sf w}_1+{\sf w}_2}\cdot (\sqrt{-1})^{{\sf w}_1}\cdot G(\omega_{\itPi_1}\omega_{\itPi_2})\cdot \Vert f_{\itPi_1}\Vert} \right)\\
&=\frac{L^{(\infty)}(m,{}^\sigma\!\itPi_1 \times {}^\sigma\!\itPi_2)}{(2\pi\sqrt{-1})^{2m+\kappa_1+{\sf w}_1+{\sf w}_2}\cdot (\sqrt{-1})^{{\sf w}_1}\cdot G({}^\sigma\!\omega_{\itPi_1}{}^\sigma\!\omega_{\itPi_2})\cdot \Vert f_{{}^\sigma\!\itPi_1}\Vert}.
\end{align*}
\item[(2)]{\rm (Garrett--Harris}\,\cite[Theorem 4.6]{GH1993}{\rm,\,C.-}\,\cite[Theorem 1.2]{Chen2021d}{\rm)}
Assume
\[
\kappa_1+\kappa_2+\kappa_3 \geq 2\max\{\kappa_1,\kappa_2,\kappa_3\}+2.
\]
For $\sigma \in {\rm Aut}(\C)$ and $m+\tfrac{1}{2} \in \Z+\tfrac{1}{2}$ critical for the triple product $L$-function $L(s,\itPi_1\times\itPi_2\times\itPi_3)$, we have
\begin{align*}
&\left( \frac{L^{(\infty)}(m+\tfrac{1}{2},\itPi_1\times\itPi_2\times\itPi_3)}{(2\pi\sqrt{-1})^{4m+\kappa_1+\kappa_2+\kappa_3+2{\sf w}'+2}\cdot(\sqrt{-1})^{{\sf w}'}\cdot G(\omega')^2\cdot \Vert f_{\itPi_1} \Vert\Vert f_{\itPi_2} \Vert\Vert f_{\itPi_3} \Vert}\right)\\
& = \frac{L^{(\infty)}(m+\tfrac{1}{2},{}^\sigma\!\itPi_1\times{}^\sigma\!\itPi_2\times{}^\sigma\!\itPi_3)}{(2\pi\sqrt{-1})^{4m+\kappa_1+\kappa_2+\kappa_3+2{\sf w}'+2}\cdot(\sqrt{-1})^{{\sf w}'}\cdot G({}^\sigma\!\omega')^2\cdot \Vert f_{{}^\sigma\!\itPi_1} \Vert\Vert f_{{}^\sigma\!\itPi_2} \Vert\Vert f_{{}^\sigma\!\itPi_3} \Vert}.
\end{align*}
Here ${\sf w}' = {\sf w}_1+{\sf w}_2+{\sf w}_3$ and $\omega' = \omega_{\itPi_1}\omega_{\itPi_2}\omega_{\itPi_3}$.
\item[(3)] Assume
\[
3\kappa_1-7 \geq \kappa_2 \geq \kappa_1+5.
\]
For $\sigma \in {\rm Aut}(\C)$ and $m \in \Z$ critical for $L(s,{\rm Sym}^3\itPi_1 \times \itPi_2)$, we have
\begin{align*}
&\sigma \left( \frac{L^{(\infty)}(m,{\rm Sym}^3\itPi_1 \times \itPi_2)}{(2\pi\sqrt{-1})^{4m+3\kappa_1+\kappa_2+6{\sf w}_1+2{\sf w}_2}\cdot(\sqrt{-1})^{{\sf w}_1+{\sf w}_2}\cdot G(\omega_{\itPi_1}^3\omega_{\itPi_2})^2\cdot \Vert f_{\itPi_1} \Vert^3\Vert f_{\itPi_2} \Vert }\right)\\
& = \frac{L^{(\infty)}(m,{\rm Sym}^3{}^\sigma\!\itPi_1 \times {}^\sigma\!\itPi_2)}{(2\pi\sqrt{-1})^{4m+3\kappa_1+\kappa_2+6{\sf w}_1+2{\sf w}_2}\cdot(\sqrt{-1})^{{\sf w}_1+{\sf w}_2}\cdot G({}^\sigma\!\omega_{\itPi_1}^3{}^\sigma\!\omega_{\itPi_2})^2\cdot \Vert f_{{}^\sigma\!\itPi_1} \Vert^3\Vert f_{{}^\sigma\!\itPi_2} \Vert }.
\end{align*}
\end{itemize}
\end{thm}

\begin{proof}
We prove the third assertion. We assume $\itPi_1$ is non-CM. The case when $\itPi_1$ is of CM-type follows from the algebraicity of critical values for $\GL_2 \times \GL_2$ due to Shimura \cite{Shimura1976}. First we recall a result of Liu \cite{Liu2019b}: Let $\itPsi$ be a cohomological irreducible cuspidal automorphic representation of $\GSp_4(\A)$ such that $\itPsi_\infty$ is a holomorphic discrete series representation of $\GSp_4(\R)$. Then the minimal ${\rm U}(2)$-type of $\itPsi_\infty$ is equal to $(k_1,k_2)$ for some $k_1 \geq k_2 \geq 3$.
Let ${\sf u}\in\Z$ be the integer such that $|\omega_\itPsi| = |\mbox{ }|_\A^{\sf u}$.
Consider the standard $L$-function $L(s,\itPsi,{\rm std})$ of $\itSigma$. Liu proved that there exists a sequence of complex numbers $(P({}^\sigma\!\itPsi))_{\sigma \in {\rm Aut}(\C)}$ such that
\begin{align}\label{T:Liu}
\begin{split}
&\sigma\left(  \frac{L^{(\infty)}(m,\itPsi,{\rm std})}{(2\pi\sqrt{-1})^{3m}\cdot (\sqrt{-1})^{\sf u}\cdot P(\itPsi)}\right) = \frac{L^{(\infty)}(m,{}^\sigma\!\itPsi,{\rm std})}{(2\pi\sqrt{-1})^{3m}\cdot (\sqrt{-1})^{\sf u}\cdot P({}^\sigma\!\itPsi)}
\end{split}
\end{align}
for all critical point $m \geq 1$ of $L(s,\itPsi,{\rm std})$ and $\sigma \in {\rm Aut}(\C)$.
On the other hand, Morimoto \cite{Morimoto2018} proved the following result on algebraicity of critical values for $\GSp_4 \times \GL_2$: Let $\itSigma$ be as above. Consider the Rankin--Selberg $L$-function $L(s,\itPsi \times \itPi_2)$.
Assume 
\[
k_1+k_2-7 \geq \kappa_2 \geq k_1-k_2+7.
\]
Then 
\begin{align}\label{T:Morimoto}
\begin{split}
&\sigma \left( \frac{L^{(\infty)}(m,\itPsi \times \itPi_2)}{(2\pi\sqrt{-1})^{4m+\kappa_2+2{\sf u}+2{\sf w}_2}\cdot (\sqrt{-1})^{{\sf u}+{\sf w}_2}\cdot G(\omega_\itPsi\omega_{\itPi_2})^2\cdot P(\itPsi)\cdot \Vert f_{\itPi_2}\Vert}\right)\\
& = \frac{L^{(\infty)}(m,{}^\sigma\!\itPsi \times {}^\sigma\!\itPi_2)}{(2\pi\sqrt{-1})^{4m+\kappa_2+2{\sf u}+2{\sf w}_2}\cdot (\sqrt{-1})^{{\sf u}+{\sf w}_2}\cdot G({}^\sigma\!\omega_\itPsi{}^\sigma\!\omega_{\itPi_2})^2\cdot P({}^\sigma\!\itPsi)\cdot \Vert f_{{}^\sigma\!\itPi_2}\Vert}
\end{split}
\end{align}
for all critical point $m \in \Z$ of $L(s,\itPsi \times \itPi_2)$ and $\sigma \in {\rm Aut}(\C)$. Note that the above algebraicity result was proved in \cite{Morimoto2018} under the assumption that ${\sf u}$ is even, and we extend the result to arbitrary ${\sf u}$ in \cite[\S\,4.5]{Chen2021f}.
Now we take $\itPsi$ be the descent of ${\rm Sym}^3 \itPi_1$ from $\GL_4(\A)$ to $\GSp_4(\A)$ with respect to the spin representation of $\GSp_4(\C)$. The existence of the descent is a consequence of Arthur's multiplicity formula for $\GSp_4(\A)$ established by Gee and Ta\"ibi \cite[Theorem 7.4.1]{GT2019}. 
Note that in this case,
\[
(k_1,k_2) = (2\kappa_1-1,\kappa_1+1),\quad {\sf u} = 3{\sf w}_1,\quad \omega_\itPsi = \omega_{\itPi_1}^3
\]
By the functoriality, we have
\[
L(s,\itPsi,{\rm std}) = L(s, {\rm Sym}^4\itPi_1\otimes \omega_{\itPi_1}^{-2}).
\]
It follows from Theorem \ref{T:Liu} and our previous result \cite{Chen2021} on Deligne's conjecture for critical values of $L(s, {\rm Sym}^4\itPi_1\otimes \omega_{\itPi_1}^{-2})$ (see also \cite{Morimoto2021}), we obtain the period relation
\[
\sigma \left( \frac{P(\itPsi)}{(2\pi\sqrt{-1})^{3\kappa_1}\cdot \Vert f_{\itPi_1}\Vert^3}\right) = \frac{P({}^\sigma\!\itPsi)}{(2\pi\sqrt{-1})^{3\kappa_1}\cdot \Vert f_{{}^\sigma\!\itPi_1}\Vert^3}
\]
for all $\sigma \in {\rm Aut}(\C)$. Combining this period relation with (\ref{T:Morimoto}), we obtain assertion (3). This completes the proof.
\end{proof}

\subsection{Proof of Theorem \ref{T:main}}

Now we begin the proof of Theorem \ref{T:main}.
Let 
\[
f(\tau) = \sum_{n=1}^\infty a_f(n)q^n \in S_\kappa(N,\omega),\quad q=e^{2\pi\sqrt{-1}\,\tau}
\]
be a normalized elliptic newform. We assume $f$ is not a CM-form.
Fix ${\sf w} \in \Z$ such that $\kappa \equiv {\sf w}\,({\rm mod}\,2)$.
Let $\itPi$ be the cohomological irreducible cuspidal automorphic representation of $\GL_2(\A)$ characterized by the following conditions:
\begin{itemize}
\item For each prime number $p \nmid N$, $\itPi_p$ is unramified and its Satake parameters are roots of the polynomial
\[
X^2-p^{(1-\kappa-{\sf w})/2}a_f(p)X+p^{-{\sf w}}\omega(p).
\]
\item $\itPi_\infty = D_\kappa \otimes |\mbox{ }|^{{\sf w}/2}$.
\end{itemize}
Let $f_\itPi \in \itPi$ be the normalized newform of $\itPi$  and $\Vert f_\itPi \Vert$ its Petersson norm in (\ref{E:Petersson norm}).
Then we have
\begin{align}\label{E:Petersson norm 2}
\Vert f_\itPi \Vert = 2\cdot\< f,f\>.
\end{align}
Recall $c^\pm(f) \in \C^\times$ is the periods of $f$ in \cite{Shimura1977}.
Let $p(\itPi,\pm)$ be the Betti-Whittaker periods of $\itPi$ with respect to the generator $[\itPi_\infty]_\pm$ following the normalization in \cite[(3.6) and (3.7)]{RT2011}.
We have the following period relation (cf.\,\cite[\S\,4.6]{RT2011}):
\begin{align}\label{E:Shimura period relation}
\sigma \left( \frac{c^\pm(f)}{p(\itPi,\pm(-1)^{(\kappa+{\sf w})/2})}\right) = \frac{c^\pm({}^\sigma\!f)}{p({}^\sigma\!\itPi,\pm(-1)^{(\kappa+{\sf w})/2})}
\end{align}
for all $\sigma \in {\rm Aut}(\C)$.
For $n \geq 1$, let ${\rm Sym}^n\itPi$ be the functorial lift of $\itPi$ to $\GL_{n+1}(\A)$ with respect to the symmetric $n$-th power representation of $\GL_2(\C)$.
The existence of the lifts was proved by Newton--Thorne \cite{NT2021} (see also \cite{GJ1972}, \cite{KS2002}, \cite{Kim2003}, \cite{CT2015}, \cite{CT2017} for $n\leq 8$).
Note that ${\rm Sym}^n\itPi$ is a cohomological irreducible cuspidal automorphic representation of $\GL_{n+1}(\A)$ (cf.\,\cite[Theorem 5.3]{Raghuram2009}).
For any primitive Dirichlet character $\chi$, we have
\begin{align}\label{E:L-function}
L^{(\infty)}(s+\tfrac{1}{2},{\rm Sym}^5\itPi \otimes\chi_\A) = L(s+\tfrac{5(\kappa+{\sf w})}{2}-2,{\rm Sym}^5(f)\otimes\chi).
\end{align}
Here $\chi_\A$ is the finite order Hecke character of $\A^\times$ such that $\chi_{\A,\,p}$ is unramified and $\chi_{\A,\,p}(p) = \chi(p)$ for all prime numbers $p$ coprime to the conductor of $\chi$.

We fix an auxiliary cohomological irreducible cuspidal automorphic representation $\itPi'$ of $\GL_2(\A)$ satisfying the following conditions:
\begin{itemize}
\item $\itPi'$ is non-CM.
\item $\itPi'_\infty = D_{2\kappa-1}\otimes |\mbox{ }|^{1/2+{\sf w}}$.
\end{itemize}
We also fix a finite order Hecke character $\eta$ of $\A^\times$ such that 
\[
\eta_\infty(-1) = (-1)^{1+{\sf w}}.
\]
Let $\itPi \boxtimes \itPi'$ be the functorial lift of the Rankin--Selberg convolution of $\itPi$ and $\itPi'$ to $\GL_4(\A)$.
The existence of the lift was proved by Ramakrishnan in \cite{Rama2000}. Since $\itPi$ and $\itPi'$ are both non-CM and the weights of $\itPi_\infty$ and $\itPi_\infty'$ are different, we see that $\itPi \boxtimes \itPi'$ is cuspidal automorphic by the cuspidality criterion in \cite[Theorem M]{Rama2000}.
Let $\itSigma_3$ and $\itSigma_4$ be tamely isobaric automorphic representations of $\GL_3(\A)$ and $\GL_4(\A)$, respectively, defined by
\begin{align*}
\itSigma_3 = (\itPi'\otimes|\mbox{ }|_\A^{-1/2}) \boxplus \eta|\mbox{ }|_\A^{\sf w},\quad \itSigma_4 = (\itPi \boxtimes \itPi')\otimes |\mbox{ }|_\A^{-1/2}.
\end{align*}
It is easy to verify that $\itSigma_3$ and $\itSigma_4$ are cohomological. Indeed, we have
\[
\itSigma_{3,\infty} = {\rm Sym}^2 \itPi_\infty,\quad \itSigma_{4,\infty} = {\rm Sym}^3 \itPi_\infty.
\]
Also note that $(\itSigma_{3,\infty},\itSigma_{4,\infty})$ is balanced (cf.\,\cite[Theorem 5.3]{Raghuram2009}).
Let $\chi$ be a finite order Hecke character of $\A^\times$.
We write
\[
\itSigma_3' = {\rm Sym}^2\itPi,\quad \itSigma_4' = {\rm Sym}^3\itPi
\]
and consider four Rankin--Selberg $L$-functions for $\GL_4 \times \GL_3$:
\begin{align}\label{E:RS L-function}
\begin{split}
&L(s,\itSigma_4 \times \itSigma_3' \otimes \chi),\quad L(s,\itSigma_4 \times \itSigma_3\otimes \chi),\\
&L(s,\itSigma_4' \times \itSigma_3'\otimes \chi),\quad L(s,\itSigma_4' \times \itSigma_3\otimes \chi).
\end{split}
\end{align}
Note that the sets of critical points for these $L$-functions are the same and is given by
\[
\left\{ m+\tfrac{1}{2} \in \Z+\tfrac{1}{2}\,\left\vert\,\tfrac{-\kappa-5{\sf w}}{2}+1 \leq m \leq \tfrac{\kappa-5{\sf w}}{2}-1\right\}\right..
\]
From now on, we assume $\kappa\geq6$.

\subsubsection{Case $m+\tfrac{1}{2} \neq \tfrac{1-5{\sf w}}{2}$}
Let $m+\tfrac{1}{2} \neq \tfrac{1-5{\sf w}}{2}$ be a non-central critical point.
In particular, since $m+\tfrac{1}{2}$ is non-central, the Rankin--Selberg $L$-functions in (\ref{E:RS L-function}) are non-vanishing at $s = m+\tfrac{1}{2}$ by the results of Jacquet--Shalika \cite[Theorem 5.3]{JS1981} and Shahidi \cite[Theorem 5.2]{Shahidi1981}. 
We have the following factorizations of $L$-functions:
\begin{align*}
L(s,\itSigma_4 \times \itSigma_3' \otimes \chi) &= L(s-\tfrac{1}{2},{\rm Sym}^3\itPi \times \itPi'\otimes \chi)\cdot L(s-\tfrac{1}{2},\itPi\times \itPi'\otimes\omega_\itPi\chi),\\
L(s,\itSigma_4 \times \itSigma_3\otimes \chi) &= L(s-1,\itPi \times \itPi' \times \itPi'\otimes \chi)\cdot L(s-\tfrac{1}{2}+{\sf w},\itPi \times \itPi'\otimes\eta\chi).
\end{align*}
By Theorem \ref{T:algebraicity results}, for $\sigma \in {\rm Aut}(\C)$ we have
\begin{align}\label{E:main proof 1}
\begin{split}
&\sigma \left( \frac{L^{(\infty)}(m+\tfrac{1}{2},\itSigma_4 \times \itSigma_3' \otimes \chi)\cdot L^{(\infty)}(m+\tfrac{1}{2},\itSigma_4 \times \itSigma_3\otimes \chi)^{-1} }{(2\pi\sqrt{-1})\cdot(\sqrt{-1})\cdot G(\omega_\itPi^3\omega_{\itPi'}^{-1}\eta^{-1})^2\cdot \Vert f_\itPi\Vert^2\Vert f_{\itPi'}\Vert^{-1}}\right)\\
& = \frac{L^{(\infty)}(m+\tfrac{1}{2},{}^\sigma\!\itSigma_4 \times {}^\sigma\!\itSigma_3' \otimes {}^\sigma\!\chi)\cdot L^{(\infty)}(m+\tfrac{1}{2},{}^\sigma\!\itSigma_4 \times {}^\sigma\!\itSigma_3\otimes {}^\sigma\!\chi)^{-1} }{(2\pi\sqrt{-1})\cdot(\sqrt{-1})\cdot G({}^\sigma\!\omega_\itPi^3{}^\sigma\!\omega_{\itPi'}^{-1}{}^\sigma\!\eta^{-1})^2\cdot \Vert f_{{}^\sigma\!\itPi}\Vert^2\Vert f_{{}^\sigma\!\itPi'}\Vert^{-1}}.
\end{split}
\end{align}
Note that the assumption $\kappa \geq 6$ is needed in order to apply Theorem \ref{T:algebraicity results}-(3).
We also have the following factorizations of $L$-functions:
\begin{align*}
L(s,\itSigma_4' \times \itSigma_3\otimes \chi) &= L(s-\tfrac{1}{2},{\rm Sym}^3\itPi \times\itPi'\otimes \chi)\cdot L(s+{\sf w},{\rm Sym}^3\itPi \otimes\eta\chi),\\
L(s,\itSigma_4' \times \itSigma_3'\otimes \chi) &= L(s,{\rm Sym}^5\itPi\otimes \chi)\cdot L(s,{\rm Sym}^3\itPi\otimes\omega_\itPi\chi)\cdot L(s, \itPi\otimes\omega_\itPi^2\chi).
\end{align*}
By \cite[Theorem 6.2]{GH1993} and \cite[Theorem 1.6]{Chen2021d} on Deligne's conjecture for symmetric cube $L$-functions, for any finite order Hecke character $\xi$ of $\A^\times$ and $\sigma \in {\rm Aut}(\C)$, we have
\begin{align*}
&\sigma \left( \frac{L^{(\infty)}(m+{\sf w}+\tfrac{1}{2},{\rm Sym}^3\itPi\otimes\xi)}{(2\pi\sqrt{-1})^{2m+2\kappa+5{\sf w}}\cdot(\sqrt{-1})^{\sf w}\cdot G(\omega_\itPi\xi)^2\cdot p(\itPi,(-1)^{m+1+{\sf w}}\xi_\infty(-1))^2\cdot \Vert f_\itPi \Vert} \right)\\
& = \frac{L^{(\infty)}(m+{\sf w}+\tfrac{1}{2},{\rm Sym}^3{}^\sigma\!\itPi\otimes{}^\sigma\!\xi)}{(2\pi\sqrt{-1})^{2m+2\kappa+5{\sf w}}\cdot(\sqrt{-1})^{\sf w}\cdot G({}^\sigma\!\omega_\itPi{}^\sigma\!\xi)^2\cdot p({}^\sigma\!\itPi,(-1)^{m+1+{\sf w}}\xi_\infty(-1))^2\cdot \Vert f_{{}^\sigma\!\itPi} \Vert}.
\end{align*}
Combining with Theorems \ref{T:GL_2} and \ref{T:algebraicity results}, for $\sigma \in {\rm Aut}(\C)$ we have
\begin{align}\label{E:main proof 2}
\begin{split}
&\sigma \left( \frac{L^{(\infty)}(m+\tfrac{1}{2},{\rm Sym}^5\itPi\otimes\chi)\cdot L^{(\infty)}(m+\tfrac{1}{2},\itSigma_4' \times \itSigma_3\otimes \chi)\cdot L^{(\infty)}(m+\tfrac{1}{2},\itSigma_4' \times \itSigma_3'\otimes \chi)^{-1}}{(2\pi\sqrt{-1})^{3m-1+(9\kappa+15{\sf w})/2}\cdot(\sqrt{-1})^{1+{\sf w}}\cdot G(\omega_{\itPi'}^2\eta^2\chi^3)\cdot p(\itPi,(-1)^m\chi_\infty(-1))^3\cdot\Vert f_\itPi\Vert^3}\right) \\
& = \frac{L^{(\infty)}(m+\tfrac{1}{2},{\rm Sym}^5{}^\sigma\!\itPi\otimes{}^\sigma\!\chi)\cdot L^{(\infty)}(m+\tfrac{1}{2},{}^\sigma\!\itSigma_4' \times {}^\sigma\!\itSigma_3\otimes {}^\sigma\!\chi)\cdot L^{(\infty)}(m+\tfrac{1}{2},{}^\sigma\!\itSigma_4' \times {}^\sigma\!\itSigma_3'\otimes {}^\sigma\!\chi)^{-1}}{(2\pi\sqrt{-1})^{3m-1+(9\kappa+15{\sf w})/2}\cdot(\sqrt{-1})^{1+{\sf w}}\cdot G({}^\sigma\!\omega_{\itPi'}^2{}^\sigma\!\eta^2{}^\sigma\!\chi^3)\cdot p({}^\sigma\!\itPi,(-1)^m\chi_\infty(-1))^3\cdot\Vert f_{{}^\sigma\!\itPi}\Vert^3}.
\end{split}
\end{align}
On the other hand, it follows from Corollary \ref{C:ratio of L-values} that
\begin{align*}
&\sigma \left( \frac{L^{(\infty)}(m+\tfrac{1}{2},\itSigma_4 \times \itSigma_3')\cdot L^{(\infty)}(m+\tfrac{1}{2},\itSigma_4' \times \itSigma_3)}{L^{(\infty)}(m+\tfrac{1}{2},\itSigma_4 \times \itSigma_3)\cdot L^{(\infty)}(m+\tfrac{1}{2},\itSigma_4' \times \itSigma_3')} \right)\\
& = \frac{L^{(\infty)}(m+\tfrac{1}{2},{}^\sigma\!\itSigma_4 \times {}^\sigma\!\itSigma_3')\cdot L^{(\infty)}(m+\tfrac{1}{2},{}^\sigma\!\itSigma_4' \times {}^\sigma\!\itSigma_3)}{L^{(\infty)}(m+\tfrac{1}{2},{}^\sigma\!\itSigma_4 \times {}^\sigma\!\itSigma_3)\cdot L^{(\infty)}(m+\tfrac{1}{2},{}^\sigma\!\itSigma_4' \times {}^\sigma\!\itSigma_3')}
\end{align*}
for all $\sigma \in {\rm Aut}(\C)$.
We thus conclude from (\ref{E:main proof 1}) and (\ref{E:main proof 2}) that 
\begin{align}\label{E:DC}
\begin{split}
&\sigma \left( \frac{L^{(\infty)}(m+\tfrac{1}{2},{\rm Sym}^5\itPi\otimes\chi)}{(2\pi\sqrt{-1})^{3m+(9\kappa+15{\sf w})/2}\cdot(\sqrt{-1})^{\sf w}\cdot G(\omega_\itPi^2\chi)^3\cdot p(\itPi,(-1)^m\chi_\infty(-1))^3\cdot \Vert f_\itPi \Vert^3} \right)\\
& = \frac{L^{(\infty)}(m+\tfrac{1}{2},{\rm Sym}^5{}^\sigma\!\itPi\otimes{}^\sigma\!\chi)}{(2\pi\sqrt{-1})^{3m+(9\kappa+15{\sf w})/2}\cdot(\sqrt{-1})^{\sf w}\cdot G({}^\sigma\!\omega_\itPi^2{}^\sigma\!\chi)^3\cdot p({}^\sigma\!\itPi,(-1)^m\chi_\infty(-1))^3\cdot \Vert f_{{}^\sigma\!\itPi} \Vert^3}
\end{split}
\end{align}
for all $\sigma \in {\rm Aut}(\C)$.
Recall the identity (\ref{E:L-function}) for $L$-functions and the period relations (\ref{E:Petersson norm 2}), (\ref{E:Shimura period relation}). 
This completes the proof of Theorem \ref{T:main} for non-central critical points.

\subsubsection{Case $m+\tfrac{1}{2} = \tfrac{1-5{\sf w}}{2}$}

Now we consider the central critical point $s=\tfrac{1-5{\sf w}}{2}$. Note that in this case ${\sf w}$ must be even. 
We have a factorization of the twisted exterior square $L$-function of ${\rm Sym}^5\itPi$:
\[
L(s,{\rm Sym}^5\itPi,\exterior{2} \otimes \,\omega_\itPi^{-5}) = L(s,{\rm Sym}^8\itPi\otimes \omega_\itPi^{-4})\cdot L(s,{\rm Sym}^4\itPi\otimes \omega_\itPi^{-2})\cdot \zeta_\F(s).
\]
In particular, $L(s,{\rm Sym}^5\itPi,\exterior{2} \otimes \,\omega_\itPi^{-5})$ has a pole at $s=1$.
In this case, for $\sigma \in {\rm Aut}(\C)$, we have the Betti--Shilika period $\omega({\rm Sym}^5 {}^\sigma\!\itPi,\pm) \in \C^\times$ of ${\rm Sym}^5 {}^\sigma\!\itPi$ in \cite[Definition/Proposition 4.2.1]{GR2014} (recalled in (\ref{A:Shalika period})). 
Fix a finite order Hecke character $\eta$ of $\A^\times$ such that $\eta_\infty(-1) = - \chi_\infty(-1)$.
By \cite[Theorem 7.1.2]{GR2014}, we have
\begin{align}\label{E:main 2 proof 1}
\begin{split}
&\sigma \left(\frac{ L^{(\infty)}(m+\tfrac{1}{2},{\rm Sym}^5\itPi \otimes \chi)\cdot L^{(\infty)}(m+\tfrac{1}{2},{\rm Sym}^5\itPi \otimes \eta)^{-1}}{ \omega({\rm Sym}^5 \itPi,(-1)^m\chi_\infty(-1))\cdot\omega({\rm Sym}^5 \itPi,(-1)^{1+m}\chi_\infty(-1))^{-1}\cdot G(\chi\eta^{-1})^3}\right)\\
& = \frac{L^{(\infty)}(m+\tfrac{1}{2},{\rm Sym}^5{}^\sigma\!\itPi \otimes {}^\sigma\!\chi)\cdot L^{(\infty)}(m+\tfrac{1}{2},{\rm Sym}^5{}^\sigma\!\itPi \otimes {}^\sigma\!\eta)^{-1}}{ \omega({\rm Sym}^5 {}^\sigma\!\itPi,(-1)^m\chi_\infty(-1))\cdot\omega({\rm Sym}^5 {}^\sigma\!\itPi,(-1)^{1+m}\chi_\infty(-1))^{-1}\cdot G({}^\sigma\!\chi{}^\sigma\!\eta^{-1})^3}
\end{split}
\end{align}
for all critical points $m+\tfrac{1}{2} \neq \tfrac{1-5{\sf w}}{2}$ and $\sigma \in {\rm Aut}(\C)$. 
On the other hand, by \cite[Theorem 7.21]{HR2020}, we have
\begin{align}\label{E:main 2 proof 2}
\begin{split}
&\sigma \left(\frac{L^{(\infty)}(m-\tfrac{1}{2},{\rm Sym}^5\itPi \otimes \chi)}{(2\pi\sqrt{-1})^{-3}\cdot\Omega({\rm Sym}^5\itPi,(-1)^{1+m}\chi_\infty(-1))\cdot L^{(\infty)}(m+\tfrac{1}{2},{\rm Sym}^5\itPi \otimes \chi)} \right)\\
& = \frac{L^{(\infty)}(m-\tfrac{1}{2},{\rm Sym}^5{}^\sigma\!\itPi \otimes {}^\sigma\!\chi)}{(2\pi\sqrt{-1})^{-3}\cdot\Omega({\rm Sym}^5{}^\sigma\!\itPi,(-1)^{1+m}\chi_\infty(-1))\cdot L^{(\infty)}(m+\tfrac{1}{2},{\rm Sym}^5{}^\sigma\!\itPi \otimes {}^\sigma\!\chi)}
\end{split}
\end{align}
for all critical points $m+\tfrac{1}{2} \neq \tfrac{1-5{\sf w}}{2}$ and $\sigma \in {\rm Aut}(\C)$, 
where $\Omega({\rm Sym}^5\itPi,\pm) \in \C^\times$ is the relative period of ${\rm Sym}^5\itPi$ defined in \cite[Definition 5.3]{HR2020} (recalled in (\ref{A:bottom relative period})). In Theorem \ref{A:period relation} below, we show that
\[
\sigma \left( \Omega({\rm Sym}^5\itPi,\pm)\cdot \frac{\omega({\rm Sym}^5\itPi,\mp)}{\omega({\rm Sym}^5\itPi,\pm)}\right) = \Omega({\rm Sym}^5{}^\sigma\!\itPi,\pm)\cdot \frac{\omega({\rm Sym}^5{}^\sigma\!\itPi,\mp)}{\omega({\rm Sym}^5{}^\sigma\!\itPi,\pm)}
\]
for all $\sigma \in {\rm Aut}(\C)$.
We then conclude from (\ref{E:main 2 proof 1}) and (\ref{E:main 2 proof 2}) that
\begin{align*}
&\sigma \left( \frac{L^{(\infty)}(\tfrac{1-5{\sf w}}{2},{\rm Sym}^5\itPi \otimes \chi)}{(2\pi\sqrt{-1})^{-3}\cdot G(\chi\eta^{-1})^3\cdot L^{(\infty)}(\tfrac{3-5{\sf w}}{2},{\rm Sym}^5\itPi \otimes \eta)}\right)\\
& = \frac{L^{(\infty)}(\tfrac{1-5{\sf w}}{2},{\rm Sym}^5{}^\sigma\!\itPi \otimes {}^\sigma\!\chi)}{(2\pi\sqrt{-1})^{-3}\cdot G({}^\sigma\!\chi{}^\sigma\!\eta^{-1})^3\cdot L^{(\infty)}(\tfrac{3-5{\sf w}}{2},{\rm Sym}^5{}^\sigma\!\itPi \otimes {}^\sigma\!\eta)}
\end{align*}
for all $\sigma \in {\rm Aut}(\C)$.
Therefore, (\ref{E:DC}) for the central critical value $L^{(\infty)}(\tfrac{1-5{\sf w}}{2},{\rm Sym}^5\itPi \otimes \chi)$ follows from (\ref{E:DC}) for $L^{(\infty)}(\tfrac{3-5{\sf w}}{2},{\rm Sym}^5\itPi \otimes \eta)$.
This completes the proof.

\begin{rmk}
In the above proof of Conjecture \ref{C:Deligne}, the argument works for $\kappa \geq 3$. 
It suffices to relax the assumption in Theorem \ref{T:algebraicity results}-(3) in order to prove Conjecture \ref{C:Deligne} for $\kappa \geq 3$.
\end{rmk}

\begin{rmk}
Since we assume $\kappa \geq 6$, $m+\tfrac{1}{2} = \frac{5-5{\sf w}}{2}$ is also a critical point when ${\sf w}$ is even. Then, by \cite[Theorem 7.21]{HR2020} alone, we can deduce (\ref{E:DC}) for the central critical value from that for $L^{(\infty)}(\tfrac{5-5{\sf w}}{2},{\rm Sym}^5\itPi \otimes \chi)$.
However, the above argument are necessary if we consider $\kappa=4$. As in this case, there are no non-central critical points with the same parity as the central critical point.
\end{rmk}

\subsection{Period relations}
We keep the notation of the previous section.
For $n \geq 1$, recall ${\rm Sym}^n\itPi$ is the symmetric $n$-th power lifting of $\itPi$. 
Note that ${\rm Sym}^n\itPi$ is a cohomological tamely isobaric automorphic representation of $\GL_{n+1}(\A)$.
Based on Deligne's conjecture for critical values of the standard $L$-function of ${\rm Sym}^n\itPi$ in \cite[Proposition 7.7]{Deligne1979}, we propose the following conjectural period relations between the periods associated to $\itPi$ and the Betti--Whittaker periods of ${\rm Sym}^n\itPi$.

\begin{conj}\label{C: period relations}
For $n \geq 1$, there exists a constant $C_{n,\infty} \in \C^\times$ unique up to $\Q^\times$, depending only on $n$ and $\itPi_\infty$, such that for all $\sigma \in {\rm Aut}(\C)$, we have
\begin{align*}
&\sigma \left( \frac{p({\rm Sym}^{2r+1}\itPi,\pm)}{C_{2r+1,\infty}\cdot G(\omega_\itPi)^{2r(r+1)^2}\cdot p(\itPi,\pm)^{r+1}\cdot \Vert f_\itPi\Vert^{2r(r+1)(2r+1)/3}} \right)\\
& = \frac{p({\rm Sym}^{2r+1}{}^\sigma\!\itPi,\pm)}{C_{2r+1,\infty}\cdot G({}^\sigma\!\omega_\itPi)^{2r(r+1)^2}\cdot p({}^\sigma\!\itPi,\pm)^{r+1}\cdot \Vert f_{{}^\sigma\!\itPi}\Vert^{2r(r+1)(2r+1)/3}}
\end{align*}
if $n=2r+1$, and
\begin{align*}
&\sigma \left( \frac{p({\rm Sym}^{2r}\itPi)}{C_{2r,\infty}\cdot G(\omega_\itPi)^{(2r+1)r^2}\cdot \Vert f_\itPi\Vert^{r(r+1)(2r+1)/3}} \right) = \frac{p({\rm Sym}^{2r}{}^\sigma\!\itPi)}{C_{2r,\infty}\cdot G({}^\sigma\!\omega_\itPi)^{(2r+1)r^2}\cdot \Vert f_{{}^\sigma\!\itPi}\Vert^{r(r+1)(2r+1)/3}}
\end{align*}
if $n=2r$.
\end{conj}

\begin{rmk}
In \cite[Theorem 4.13]{Chen2020}, we prove the conjecture for $n=2$ and explicitly determine the constant $C_{2,\infty}$.
\end{rmk}

As a consequence of Theorem \ref{T:main}, we prove the conjecture for $n=3$.
\begin{thm}\label{A: period relations 2}
If $\kappa \geq 6$, then Conjecture \ref{C: period relations} holds for $n=3$.
\end{thm}

\begin{rmk}
It is customary in the literature to choose ${\sf w}=-\kappa$. In this case, we let $p({\rm Sym}^3(f),\pm)=p({\rm Sym}^3\itPi,\pm)$. Then Theorem \ref{A: period relations 2} is equivalent to Theorem \ref{T:main 2}.
\end{rmk}

\begin{proof}
Let $\chi$ be any finite order Hecke character of $\A^\times$.
We have the following factorizations of $L$-functions:
\begin{align*}
L(s,{\rm Sym}^3\itPi \times {\rm Sym}^2\itPi\otimes\chi) &= L(s,{\rm Sym}^5\itPi\otimes \chi)\cdot L(s,{\rm Sym}^3\itPi\otimes\omega_\itPi\chi)\cdot L(s, \itPi\otimes\omega_\itPi^2\chi),\\
L(s,{\rm Sym}^2\itPi \times \itPi\otimes\omega_\itPi\chi) &= L(s,{\rm Sym}^3\itPi\otimes\omega_\itPi\chi)\cdot L(s, \itPi\otimes\omega_\itPi^2\chi).
\end{align*}
Let $m+\tfrac{1}{2}$ be any non-central critical point for ${\rm Sym}^5\itPi$. By Theorem \ref{T:main}, for all $\sigma \in {\rm Aut}(\C)$ we have
\begin{align*}
&\sigma\left( \frac{L^{(\infty)}(m+\tfrac{1}{2},{\rm Sym}^3\itPi \times {\rm Sym}^2\itPi\otimes\chi)\cdot L^{(\infty)}(m+\tfrac{1}{2},{\rm Sym}^2\itPi \times \itPi\otimes\omega_\itPi\chi)^{-1}}{(2\pi\sqrt{-1})^{3m+(9\kappa+15{\sf w})/2}\cdot(\sqrt{-1})^{\sf w}\cdot G(\omega_\itPi^2\chi)^3\cdot p(\itPi,(-1)^m\chi_\infty(-1))^3\cdot \Vert f_\itPi \Vert^3} \right)\\
& = \frac{L^{(\infty)}(m+\tfrac{1}{2},{\rm Sym}^3{}^\sigma\!\itPi \times {\rm Sym}^2{}^\sigma\!\itPi\otimes{}^\sigma\!\chi)\cdot L^{(\infty)}(m+\tfrac{1}{2},{\rm Sym}^2{}^\sigma\!\itPi \times {}^\sigma\!\itPi\otimes{}^\sigma\!\omega_\itPi{}^\sigma\!\chi)^{-1}}{(2\pi\sqrt{-1})^{3m+(9\kappa+15{\sf w})/2}\cdot(\sqrt{-1})^{\sf w}\cdot G({}^\sigma\!\omega_\itPi^2{}^\sigma\!\chi)^3\cdot p({}^\sigma\!\itPi,(-1)^m\chi_\infty(-1))^3\cdot \Vert f_{{}^\sigma\!\itPi} \Vert^3}.
\end{align*}
Let
\[
C_{3,\infty} = (2\pi\sqrt{-1})^{3m+1+(9\kappa+15{\sf w})/2}\cdot \frac{Z(m,{\rm Sym}^3\itPi_\infty,{\rm Sym}^2\itPi_\infty)}{Z(m+{\sf w},{\rm Sym}^2\itPi_\infty,\itPi_\infty)}.
\]
By Theorem \ref{T:Raghuram}, it is easy to verify that Conjecture \ref{C: period relations} holds for this choice of $C_{3,\infty}$.
\end{proof}

\appendix
\section{Period relations between relative periods and Betti--Shalika periods}\label{S:appendix}

Let $\itPi$ be a cohomological irreducible cuspidal automorphic representation of $\GL_n(\A)$. We assume $n$ is even and fix a non-trivial representative $\delta$ of $K_n / K_n^\circ$.
In the appendix, we establish period relations between relative periods and ratios of Betti--Shalika periods of $\itPi$.
First we recall the definition of the relative periods $\Omega(\itPi,\pm)$ of $\itPi$ (cf.\,\cite[\S\,5.2]{HR2020}). Recall $b_{n} = \tfrac{n^2}{4}$. We put $t_{n} = \tfrac{n^2}{4}+\tfrac{n}{2}-1$. Let $M$ the irreducible algebraic representation of $\GL_{n}$ such that
\[
H^{b_n}(\frak{g}_n,K_n^\circ;\itPi\otimes M_\C) \neq 0,\quad H^{t_n}(\frak{g}_n,K_n^\circ;\itPi\otimes M_\C) \neq 0.
\]
Since $n$ is even, the above cohomology groups are $2$-dimensional. Moreover, the restriction to $(\frak{g}_n,K_n^\circ)$ of $\itPi_\infty$ decomposes as $\itPi_\infty = \itPi_\infty^{+}\oplus \itPi_\infty^{-}$ such that $\delta$ interchanges the two summands.
We fix generators
\[
c_b^{++} \in H^{b_n}(\frak{g}_n,K_n^\circ;\itPi_\infty^+\otimes M_\C),\quad c_t^{++} \in H^{t_n}(\frak{g}_n,K_n^\circ;\itPi_\infty^+\otimes M_\C).
\]
We then define generators
\[
[\itPi_\infty]_\pm^b \in H^{b_n}(\frak{g}_n,K_n^\circ;\itPi_\infty\otimes M_\C)[\pm],\quad [\itPi_\infty]_\pm^t \in H^{t_n}(\frak{g}_n,K_n^\circ;\itPi_\infty\otimes M_\C)[\pm]
\]
by
\begin{align}\label{A:generators}
[\itPi_\infty]_\pm^b = c_b^{++} \pm \delta\cdot c_b^{++},\quad [\itPi_\infty]_\pm^t = c_t^{++} \pm \delta\cdot c_t^{++}.
\end{align}
Let 
\[
T_\itPi^\pm : H^{b_n}(\frak{g}_n,K_n^\circ;\itPi_\infty\otimes M_\C)[\pm] \longrightarrow  H^{b_n}(\frak{g}_n,K_n^\circ;\itPi_\infty\otimes M_\C)[\mp]
\]
be the $\GL_n(\A_f)$-equivariant isomorphism defined as follows: Let $\Phi : \itPi_\infty \otimes \itPi_f \rightarrow \itPi$ be a $((\frak{g}_n,K_n)\times\GL_n(\A_f))$-equivariant isomorphism. It induces a $(K_n\times\GL_n(\A_f))$-isomorphism
\[
\Phi : H^{b_n}(\frak{g}_n,K_n^\circ;\itPi_\infty\otimes M_\C) \otimes \itPi_f \longrightarrow H^{b_n}(\frak{g}_n,K_n^\circ;\itPi\otimes M_\C).
\]
Then $T_\itPi^\pm$ is defined so that
\[
T_\itPi^\pm\circ\Phi([\itPi_\infty]_\pm^b\otimes{\bf v}) = \Phi([\itPi_\infty]_\mp^b\otimes{\bf v}),\quad {\bf v}\in\itPi_f.
\]
It is clear that $T_\itPi^\pm$ is independent of the choice of $\Phi$.
The relative period $\Omega(\itPi,\pm)$ is the non-zero complex number, unique up to $\Q(\itPi)^\times$, such that 
\begin{align}\label{A:bottom relative period}
T_\itPi^\pm \left(H^{b_n}(\frak{g}_n,K_n^\circ;\itPi\otimes M_\C)[\pm]^{{\rm Aut}(\C/\Q(\itPi))}\right) = \frac{H^{b_n}(\frak{g}_n,K_n^\circ;\itPi\otimes M_\C)[\mp]^{{\rm Aut}(\C/\Q(\itPi))}}{\Omega(\itPi,\pm)}.
\end{align}
For $\sigma \in {\rm Aut}(\C)$, similarly we can define $T_{{}^\sigma\!\itPi}^\pm$ and $\Omega({}^\sigma\!\itPi,\pm)$. We normalize the relative periods so that the diagram
\[
\begin{tikzcd}
H^{b_n}(\frak{g}_n,K_n^\circ;\itPi\otimes M_\C)[\pm] \arrow[rr, "{\Omega(\itPi,\pm)}\cdot T_\itPi^\pm"] \arrow[d, "\sigma^*"] & &H^{b_n}(\frak{g}_n,K_n^\circ;\itPi\otimes M_\C)[\mp]\arrow[d, "\sigma^*"]\\
H^{b_n}(\frak{g}_n,K_n^\circ;{}^\sigma\!\itPi\otimes M_\C)[\pm] \arrow[rr, "{\Omega({}^\sigma\!\itPi,\pm)}\cdot T_{{}^\sigma\!\itPi}^\pm"] & & H^{b_n}\left(\frak{g}_n,K_n^\circ;{}^\sigma\!\itPi \otimes M_\C\right)[\mp]
\end{tikzcd}
\]
commutes.
In other words, we have
\[
\sigma^*\circ T_{\itPi}^\pm \left( \Omega(\itPi,\pm)\cdot c \right) = T_{{}^\sigma\!\itPi}^\pm\left( \Omega({}^\sigma\!\itPi,\pm)\cdot \sigma^*c \right)
\]
for all $c \in H^{b_n}(\frak{g}_n,K_n^\circ;\itPi\otimes M_\C)[\pm]$.
Note that in the notation and normalization of \cite[\S\,5.2.3, Definition 5.3]{HR2020}, we have
\[
\Omega({}^\sigma\!\itPi,\pm) = (\sqrt{-1})^{n/2}\cdot \Omega^\pm({}^\sigma\!\lambda,{}^\sigma\!\itPi_f).
\]
Now we recall the definition of Betti--Shalika periods $\omega(\itPi,\pm)$ of $\itPi$.
Assume there exists an algebraic Hecke character $\eta$ of $\A^\times$ such that the twisted exterior square $L$-function 
\[
L(s,\itPi,\exterior{2}\otimes\,\eta^{-1})
\]
has a pole at $s=1$. By the result of Jacquet--Shalika \cite{JS1990b}, the assumption is equivalent to saying that the $(\eta,\psi)$-Shalika functional is non-vanishing on $\itPi$. 
For $\varphi \in \itPi$, let $S_\varphi$ be the $(\eta,\psi)$-Shalika function of $\varphi$ (cf.\,\cite[\S\,3.1]{GR2014}).
Let $\mathcal{S}(\itPi_\infty)$ and $\mathcal{S}(\itPi_f)$ be the spaces of $(\eta_\infty,\psi_\infty)$-Shalika functions and $(\eta_f,\psi_f)$-Shalika functions of $\itPi_\infty$ and $\itPi_f$, respectively. For $S_\infty \in \mathcal{S}(\itPi_\infty)$ and $S_f \in \mathcal{S}(\itPi_f)$, there exists a unique cuap form $\varphi \in \itPi$ such that
\[
S_\varphi = S_\infty\cdot S_f.
\]
In this way, we obtain a $((\frak{g}_n,K_n)\times \GL_n(\A_f))$-equivariant isomorphism
\begin{align}\label{E:Shalika isomorphism}
\mathcal{S}(\itPi_\infty) \otimes \mathcal{S}(\itPi_f) \longrightarrow \itPi.
\end{align}
Let 
\begin{align}\label{A:Shalika pm}
\Phi_{\itPi,\pm}^{\mathcal{S}} : \mathcal{S}(\itPi_f) \longrightarrow H^{t_n}(\frak{g}_n,K_n^\circ;\itPi \otimes M_\C)[\pm]
\end{align}
be the $\GL_n(\A_f)$-equivariant isomorphism defined as follows: For $S \in \mathcal{S}(\itPi_f)$, we have
\[
[\itPi_\infty]_\pm^t\otimes S \in H^{t_n}(\frak{g}_n,K_n^\circ; \mathcal{S}(\itPi_\infty) \otimes \mathcal{S}(\itPi_f)\otimes M_\C)[\pm].
\] 
Then $\Phi_{\itPi,\pm}^{\mathcal{S}}(S)$ is the image of $[\itPi_\infty]_\pm^t\otimes S$ under the $\GL_n(\A_f)$-equivariant isomorphism induced by the isomorphism (\ref{E:Shalika isomorphism}).
Here $[\itPi_\infty]_\pm^t$ is the generator fixed in (\ref{A:generators}).
For $\sigma \in {\rm Aut}(\C)$, let $t_\sigma : \mathcal{S}(\itPi_f) \rightarrow \mathcal{S}({}^\sigma\!\itPi_f)$ be the $\sigma$-linear $\GL_n(\A_f)$-equivariant isomorphism defined by
\[
t_\sigma S(g) = \sigma \left( S({\rm diag}(u_\sigma^{-1},\cdots,u_\sigma^{-1},1,\cdots,1)g)\right),\quad g \in \GL_n(\A_f).
\]
Here $1$ appears in the diagonal matrix ${n}/{2}$-times and $u_\sigma \in \widehat{\Z}^\times$ is the unique element such that $\sigma(\psi(x)) = \psi(u_\sigma x)$ for all $x \in \A_f$.
The Betti--Shalika period $\omega(\itPi,\pm)$ is the non-zero complex number, unique up to $\Q(\itPi,\eta)^\times$, such that
\[
\Phi_{\itPi,\pm}^{\mathcal{S}}\left(\frac{\mathcal{S}(\itPi_f)^{{\rm Aut}(\C/\Q(\itPi,\eta))}}{\omega(\itPi,\pm)}\right) = H^{t_n}(\frak{g}_n,K_n^\circ;\itPi \otimes M_\C)[\pm]^{{\rm Aut}(\C/\Q(\itPi,\eta))}.
\]
For $\sigma \in {\rm Aut}(\C)$, similarly we can define $\Phi_{{}^\sigma\!\itPi,\pm}^{\mathcal{S}}$ and $\omega({}^\sigma\!\itPi,\pm)$. We normalize the Betti--Shalika periods such that
\begin{align}\label{A:Shalika period}
\sigma^* \left(\frac{\Phi_{\itPi,\pm}^{\mathcal{S}}(S)}{\omega(\itPi,\pm)}\right) = \frac{\Phi_{{}^\sigma\!\itPi,\pm}^{\mathcal{S}}(t_\sigma S)}{\omega({}^\sigma\!\itPi,\pm)}
\end{align}
for all $S \in \mathcal{S}(\itPi_f)$.
As the main result of the appendix, we have the following period relation.
\begin{thm}\label{A:period relation}
For $\sigma \in {\rm Aut}(\C)$, we have
\[
\sigma \left( \Omega(\itPi,\pm)\cdot \frac{\omega(\itPi,\mp)}{\omega(\itPi,\pm)}\right) = \Omega({}^\sigma\!\itPi,\pm)\cdot \frac{\omega({}^\sigma\!\itPi,\mp)}{\omega({}^\sigma\!\itPi,\pm)}.
\]
\end{thm}

\begin{proof}
Let 
\begin{align}\label{A:Whittaker pm}
\Phi_{\itPi,\pm}^{\mathcal{W}} : \mathcal{W}(\itPi_f) \longrightarrow H^{t_n}(\frak{g}_n,K_n^\circ;\itPi \otimes M_\C)[\pm]
\end{align}
be the $\GL_n(\A_f)$-equivariant isomorphism defined as above with respect to the generator $[\itPi_\infty]_\pm^t$ and the isomorphism (\ref{E:Whittaker isomorphism}).
Similarly we define the (top degree) Betti--Whittaker period $q(\itPi,\pm)$ of $\itPi$ such that
\begin{align}\label{A:Whittaker period}
\sigma^* \left(\frac{\Phi_{\itPi,\pm}^{\mathcal{W}}(W)}{q(\itPi,\pm)}\right) = \frac{\Phi_{{}^\sigma\!\itPi,\pm}^{\mathcal{W}}(t_\sigma W)}{q({}^\sigma\!\itPi,\pm)}
\end{align}
for all $W \in \mathcal{W}(\itPi_f)$.
By the period relations in \cite[(4.6)]{RS2008} and \cite[Corollary 3.3.14]{BR2017}, we have
\[
\sigma \left( \frac{p(\itPi,\pm)\cdot q(\itPi,\mp)}{p(\itPi,\mp)\cdot q(\itPi,\pm)} \right) = \frac{p({}^\sigma\!\itPi,\pm)\cdot q({}^\sigma\!\itPi,\mp)}{p({}^\sigma\!\itPi,\mp)\cdot q({}^\sigma\!\itPi,\pm)} 
\]
for all $\sigma \in {\rm Aut}(\C)$.
On the other hand, it is clear from the definition and normalization of the relative periods and (bottom degree) Betti--Whittaker periods in \S\,\ref{SS:Whittaker periods} that 
\[
\sigma \left( \Omega(\itPi,\pm)\cdot \frac{p(\itPi,\mp)}{p(\itPi,\pm)}\right) = \Omega({}^\sigma\!\itPi,\pm)\cdot \frac{p({}^\sigma\!\itPi,\mp)}{p({}^\sigma\!\itPi,\pm)}
\]
for all $\sigma \in {\rm Aut}(\C)$ (see the explanation in \cite[\S\,5.2.3]{HR2020}). 
Therefore, to prove the assertion, it suffices to establish the following period relation:
\begin{align}\label{A:main period relation}
\sigma \left( \frac{q(\itPi,\pm)\cdot \omega(\itPi,\mp)}{q(\itPi,\mp)\cdot \omega(\itPi,\pm)} \right) = \frac{q({}^\sigma\!\itPi,\pm)\cdot \omega({}^\sigma\!\itPi,\mp)}{q({}^\sigma\!\itPi,\mp)\cdot \omega({}^\sigma\!\itPi,\pm)} 
\end{align}
for all $\sigma \in {\rm Aut}(\C)$.
Let $\itPi^+$ be the $((\frak{g}_n,K_n^\circ)\times\GL_n(\A_f))$-submodule of $\itPi$ consisting of cusp forms in $\itPi$ with archimedean components in $\itPi_\infty^+$.
Let
\begin{align*}
&\iota_{\itPi}^{\mathcal{S}} : \mathcal{S}(\itPi_f) \longrightarrow H^{t_n}(\frak{g}_n,K_n^\circ;\itPi^+ \otimes M_\C),\quad\iota_{\itPi}^{\mathcal{W}} : \mathcal{W}(\itPi_f) \longrightarrow H^{t_n}(\frak{g}_n,K_n^\circ;\itPi^+ \otimes M_\C)
\end{align*}
be the $\GL_n(\A_f)$-equivariant isomorphisms defined similar to (\ref{A:Shalika pm}) and (\ref{A:Whittaker pm}) with $[\itPi_\infty]_\pm^t$ and $H^{t_n}(\frak{g}_n,K_n^\circ;\itPi \otimes M_\C)[\pm]$ replaced by $c_t^{++}$ and $H^{t_n}(\frak{g}_n,K_n^\circ;\itPi^+ \otimes M_\C)$, respectively.
Let $C_{\mathcal{W}}^{\mathcal{S}}(\itPi)$ be the non-zero complex number, unique up to $\Q(\itPi,\eta)^\times$, such that
\[
\frac{\iota_\itPi^{\mathcal{W}}\left(\mathcal{W}(\itPi_f)^{{\rm Aut}(\C/\Q(\itPi,\eta))}\right)}{C_{\mathcal{W}}^{\mathcal{S}}(\itPi)} = \iota_\itPi^{\mathcal{S}}\left(\mathcal{S}(\itPi_f)^{{\rm Aut}(\C/\Q(\itPi,\eta))}\right).
\]
For $\sigma \in {\rm Aut}(\C)$, similarly we define $C_{\mathcal{W}}^{\mathcal{S}}({}^\sigma\!\itPi)$. We normalize these periods so that
\begin{align}\label{A:relative period}
t_\sigma \left( \frac{(\iota_\itPi^{\mathcal{S}})^{-1}\circ\iota_\itPi^{\mathcal{W}}(W)}{C_{\mathcal{W}}^{\mathcal{S}}(\itPi)}\right) = \frac{(\iota_{{}^\sigma\!\itPi}^{\mathcal{S}})^{-1}\circ\iota_{{}^\sigma\!\itPi}^{\mathcal{W}}(t_\sigma W)}{C_{\mathcal{W}}^{\mathcal{S}}({}^\sigma\!\itPi)}
\end{align}
for all $W \in \mathcal{W}(\itPi_f)$.
On the other hand, by definition we have
\[
\Phi_{\itPi,\pm}^{\mathcal{S}} = \iota_\itPi^{\mathcal{S}} \pm \delta\cdot \iota_\itPi^{\mathcal{S}},\quad \Phi_{\itPi,\pm}^{\mathcal{W}} = \iota_\itPi^{\mathcal{W}} \pm \delta\cdot \iota_\itPi^{\mathcal{W}}.
\]
Therefore, for $\sigma \in {\rm Aut}(\C)$ and $W \in \mathcal{W}(\itPi_f)$, we have
\begin{align*}
\sigma^*\left(\frac{\Phi_{\itPi,\pm}^{\mathcal{W}}( W)}{C_{\mathcal{W}}^{\mathcal{S}}(\itPi)\cdot\omega(\itPi,\pm)}\right) & = \sigma^*\left(\frac{\Phi_{\itPi,\pm}^{\mathcal{S}}\left((\iota_{\itPi}^{\mathcal{S}})^{-1}\circ\iota_{\itPi}^{\mathcal{W}}(W)\right)}{C_{\mathcal{W}}^{\mathcal{S}}(\itPi)\cdot\omega(\itPi,\pm)}\right)\\
& = \Phi_{{}^\sigma\!\itPi,\pm}^{\mathcal{S}}\left(t_\sigma \left(\frac{(\iota_\itPi^{\mathcal{S}})^{-1}\circ\iota_\itPi^{\mathcal{W}}(W)}{C_{\mathcal{W}}^{\mathcal{S}}(\itPi)}\right)\right)\cdot \frac{1}{\omega({}^\sigma\!\itPi,\pm)}\\
& = \frac{\Phi_{{}^\sigma\!\itPi,\pm}^{\mathcal{S}}\left((\iota_{{}^\sigma\!\itPi}^{\mathcal{S}})^{-1}\circ\iota_{{}^\sigma\!\itPi}^{\mathcal{W}}(t_\sigma W)\right)}{C_{\mathcal{W}}^{\mathcal{S}}({}^\sigma\!\itPi)\cdot\omega({}^\sigma\!\itPi,\pm)}\\
& = \frac{\Phi_{{}^\sigma\!\itPi,\pm}^{\mathcal{W}}(t_\sigma W)}{C_{\mathcal{W}}^{\mathcal{S}}({}^\sigma\!\itPi)\cdot\omega({}^\sigma\!\itPi,\pm)}.
\end{align*}
Here the second and third equalities follow from  (\ref{A:Shalika period}) and (\ref{A:relative period}), respectively.
Comparing with (\ref{A:Whittaker period}), we deduce that
\[
\sigma\left(\frac{q(\itPi,\pm)}{C_{\mathcal{W}}^{\mathcal{S}}(\itPi)\cdot \omega(\itPi,\pm)}\right) = \frac{q({}^\sigma\!\itPi,\pm)}{C_{\mathcal{W}}^{\mathcal{S}}({}^\sigma\!\itPi)\cdot \omega({}^\sigma\!\itPi,\pm)}
\]
for all $\sigma \in {\rm Aut}(\C)$. In particular, (\ref{A:main period relation}) follows immediately.
This completes the proof.
\end{proof}


\end{document}